\input amstex
\documentstyle {amsppt}
\magnification=1200
\vsize=9.5truein
\hsize=6.5truein
\nopagenumbers
\nologo

\def\norm#1{\left\Vert #1 \right\Vert}
\def\abs#1{\left\vert #1 \right\vert}

%Copie du 26/8 2004

\topmatter

\centerline{
International Conference on Group Theory -- 1--6 June 2003 -- Gaeta
}

\bigskip

\bigskip

\title
Cubature formulas, geometrical designs, 
\\
reproducing kernels, and Markov operators
\endtitle

\leftheadtext{Pierre de la Harpe and Claude Pache}
\rightheadtext{Cubature formulas}

\author
Pierre de la Harpe and Claude Pache
\endauthor

\thanks
The authors acknowledge support from the 
{\it Swiss National Science Foundation}.
\endthanks

\abstract
Cubature formulas and geometrical designs are described
in terms of reproducing kernels for Hilbert spaces of functions
on the one hand, 
and Markov operators associated to orthogonal group representations
on the other hand.
In this way,
several known results for spheres in Euclidean spaces,
involving cubature formulas for polynomial functions 
and spherical designs, 
are shown to generalize 
to large classes of finite measure spaces $(\Omega,\sigma)$ 
and appropriate spaces of functions inside $L^2(\Omega,\sigma)$.
The last section points out how spherical designs are related 
to a class of reflection groups 
which are (in general dense) subgroups of orthogonal groups.
\endabstract


%%{22D10, 46E22, 43A65}
\endsubjclass

\address
Pierre de la Harpe and Claude Pache, Section de Math\'ematiques, 
Universit\'e de Gen\`eve, C.P. 240, 
CH-1211 Gen\`eve 24, Suisse. \newline 
E-mail: Pierre.delaHarpe\@math.unige.ch, Claude.Pache\@math.unige.ch
\endaddress

\endtopmatter

\document

   Let $(\Omega,\sigma)$ be a finite measure space 
and let $\Cal F$ be a vector space of 
integrable real-valued functions on $\Omega$. 
It is a natural question to ask when and how integrals
$\int_{\Omega} \varphi (\omega) d\sigma(\omega)$ can be computed, or
approximated, by sums $\sum_{x\in X} W(x) \varphi (x)$, where $X$ is a
subset of $\Omega$ and $W : X \longrightarrow \Bbb R_+^*$ a weight
function, for all $\varphi \in \Cal F$.

When $\Omega$ is an interval of the real line, 
this is a basic problem of numerical integration with a glorious 
list of contributors:
Newton (1671), Cotes (1711), Simpson (1743), Gauss (1814), 
Chebyshev (1874),  Christoffel (1877), and Bernstein (1937), 
to quote but a few; 
see \cite{Gauts--97} for some historical notes,
\cite{Gauts--76}, and \cite{DavRa--84}.
The theory is inseparable of that of orthogonal polynomials
\cite{Szeg\"o--39}.

When $\Omega$ is a space of larger dimension, 
the problems involved are often of 
geometrical and combinatorial interest.
One important case is that of spheres in Euclidean spaces, 
with rotation-invariant proba\-bility measure,
as in \cite{DeGoS--77} and \cite{GoeSe--79}.
Work related to integrations domains with
$\operatorname{dim}(\Omega) > 1$ goes back to
\cite{Maxwe--77} and \cite{Appel--90};
see also the result of Voronoi (1908)
recalled in Item~1.15 below. 
Examples which have been considered include
various domains in Euclidean spaces 
(hypercubes, simplices, ..., Item~1.20),
surfaces (e.g. tori),
and Euclidean spaces with Gaussian measures (Item~3.11).

There are also interesting cases
where $\Omega$ itself is a finite set \cite{Delsa--78}. 

\medskip

   Section 1 collects the relevant definitions 
for the general case $(\Omega,\sigma)$. 
It reviews several known examples 
on intervals and spheres.

   Our main point is to show 
that there are two notions which
are convenient for the study of cubature formulas, 
even if they are rarely explicitely introduced 
in papers of this subject.  

   First, we introduce in Section 2 the formalism of
{\it reproducing kernel Hilbert spaces} 
(\cite{Arons--50}, \cite{Krein--63}, \cite{BeHaG--03}), 
which is appropriate for generalizing to other spaces 
various results which are standard for spheres. 
See for example the existence theorem for cubature formulas $(X,W)$
with bound on $\abs{X}$ (Items~2.6--2.8),
Propositions~2.10  \&~3.4 on tight cubature formulas being
geometrical  designs, 
Proposition~2.12 on the set of distances
$
D_X \, = \, \left\{ c \in \mathopen]0,\infty\mathclose[ 
\ \big\vert  \
c = d(x,y) \  \text{for some} \ x,y \in X, x \ne y \right\}
$
if $\Omega$ is a metric space 
and if \lq\lq Condition (M)\rq\rq \ is satisfied,
and Items~2.12--2.13 for a single equality which guarantees
that a pair $(X,W)$ is a cubature formula of strength $k$.

   Section 3 concentrates on spheres,  
on consequences of the existence of the central symmetry,
and in particular on spherical designs of odd strengths.
An example of \lq\lq lattice construction\rq\rq \ (Item~3.10)
provides a cubature formula in $\Bbb S^7$ of strength $11$
and size $2400$.
Item~3.11 hints at some 
very rudimentary constructions of {\it Gaussian designs}.

   Secondly, we consider in Section 4 {\it Markov operators}
(called elsewhere  
\lq\lq convolution ope\-rators\rq\rq, 
\lq\lq difference operators\rq\rq,
\lq\lq Hecke operators\rq\rq , or
\lq\lq discrete Laplace operators\rq\rq). 
On the one hand, they provide an alternative definition of cubature
formulas on spheres \cite{Pache--04}, 
and at least part of this can be generalized
to other Riemannian symme\-tric spaces of the compact type
(work in progress);
also, they show a natural connection with the work of
Lubotzky, Phillips, and Sarnak (see Theorem~4.7).
On the other hand, they suggest an interesting class of examples of
infinite {\it reflection groups,} 
as shown in Section~5 where we state and prove
an unpublished result of B.\ Venkov on remarkable sets of generators of
the orthogonal group $O(8,\Bbb Z[1/2])$.

   The subject of cubature formulas has natural connections with
several other subjects. We give brief mentions to:
\roster
\item"$\circ$" representations of finite groups 
   (see \cite{Banna79}, \cite{HarPa--04}, and Example~1.14),
\item"$\circ$" lattices in Euclidean spaces
   (see \cite{Venko--84}, \cite{VenMa--01},  and Examples~1.16 
    \& 3.10), 
\item"$\circ$" Lehmer's conjecture (Example~1.16) and
   Waring's problem (Item~3.9) from number theory, 
\item"$\circ$" Dvoretzky's theorem from Banach space theory 
   (see \cite{LyuVa--93} and Proposition~3.12). 
\endroster
\medskip

   It is a pleasure to thank Boris Venkov,
as well as B\'ela Bajnok, Eiichi Bannai, Ga\"el Collinet, and Martin Gander
for many conversations and indications
which have been most helpful during this work.

\bigskip
\head
{\bf 
1.\ Cubature formulas, designs, and polynomial spaces
}
\endhead
\medskip

   Let $(\Omega,\sigma)$ be a finite measure space. 
By a {\bf space of functions} 
in $L^2(\Omega,\sigma)$, 
we mean a linear subspace of genuine 
\footnote{
Not functions modulo equality almost everywhere!
}
real-valued functions on $\Omega$ given with an
embedding in $L^2(\Omega,\sigma)$.

\bigskip

\proclaim{1.1.\ Definition} A {\bf cubature formula} for a finite
dimensional space $\Cal F$ of functions in $L^2(\Omega,\sigma)$ is
a pair $(X,W)$, where $X$ is a finite subset of $\Omega$ and 
$W : X \longrightarrow \Bbb R_+^*$ is a weight function, such that
$$
\sum_{x\in X} W(x) \varphi (x) \, = \,
\int_{\Omega} \varphi (\omega) d\sigma(\omega)
$$
for all $\varphi \in \Cal F$.
If $W : x \longmapsto \sigma(\Omega)/\abs{X}$ 
is the uniform weight, 
the set $X$ is called a {\bf design} for $\Cal F$.
\endproclaim

   This definition would make sense for
$\Cal F \subset L^1(\Omega,\sigma)$.
But the examples we have in mind are in $L^2(\Omega,\sigma)$,
and we will use Hilbert space methods for which
we have to assume that
$\Cal F \subset L^2(\Omega,\sigma)$;
two functions $f_1,f_2 \in \Cal F$ have a scalar product
given by
$$
\langle f_1 \mid f_2 \rangle \, = \, 
\int_{\Omega} f_1(\omega)f_2(\omega)d\sigma(\omega) .
$$
There are other denominations for these and related notions, 
including \lq\lq Chebyshev qua\-drature rule\rq\rq \ 
(for our \lq\lq designs\rq\rq ) 
with appropriate \lq\lq algebraic degree of exactness\rq\rq ;  
see, e.g., \cite{Gauts--76}.

We write $\langle \cdot \mid \cdot \rangle$
to denote both the scalar product of two vectors in a Hilbert space of
functions and the scalar product in the standard Euclidean space $\Bbb R^n$.

\bigskip

\noindent {\bf 1.2.\ Example (Simpson's formula).} 
In the case of an interval
$\Omega = [a,b]$ of the real line ($a < b$)
and Lebesgue measure,  
{\it Simpson's formula}
$$
\frac{b-a}{6}\varphi \big(a\big) + 
\frac{4(b-a)}{6}\varphi \big(\frac{a+b}{2}\big) +
\frac{b-a}{6}\varphi (b) \, = \,
\int_{a}^b \varphi (t) dt
$$
provides a cubature formula for the space $\Cal F^{(3)}(a,b)$ of
polynomial functions of degree at most $3$ on $[a,b]$.
More generally, if $\varphi$ is of class $\Cal C^{4}$ on $[a,b]$,
we have
$$
\abs{
\frac{b-a}{6}\varphi \big(a\big) + 
\frac{4(b-a)}{6}\varphi \big(\frac{a+b}{2}\big) +
\frac{b-a}{6}\varphi (b) \, - \,
\int_{a}^b \varphi (t) dt
}
\, \le \,
\frac{(b-a)^5}{2880} \sup_{a \le t \le b} 
\abs{ \varphi^{(4)}(t) } .
$$
Simpson's formula is most often used in its {\it compound form},
namely on the $n$ subintervals
$[a+\frac{j-1}{n}(b-a),a+\frac{j}{n}(b-a)]$, $1 \le j \le n$,
for $n$ large enough;
in this form, it provides a very efficient tool in numerical analysis.

\bigskip

   The following notion borrows some of the ingredients of the
notion of \lq\lq polynomial space\rq\rq \ from \cite{Godsi--93}.

\bigskip

\proclaim{1.3.\ Definition} A {\bf sequence of polynomial spaces} on
a finite measure space $(\Omega,\sigma)$ is a nested sequence of
finite dimensional spaces
$$
\Bbb R = \Cal F^{(0)} \, \subset \, 
\Cal F^{(1)} \, \subset \, \cdots \, \subset \, 
\Cal F^{(k)} \, \subset \, \cdots
$$
of functions in $L^2(\Omega,\sigma)$ with the following property:
$\Cal F^{(k)}$ is linearly generated by 
products $\varphi_1\varphi_2$ 
with $\varphi_1 \in \Cal F^{(1)}$
and $\varphi_2 \in \Cal F^{(k-1)}$, 
for all $k \ge 1$.
\endproclaim

\bigskip

\noindent
{\bf 1.4.\ Remarks.} (i) For $(\Cal F^{(k)})_{k \ge 0}$ as in the
definition, there is a natural mapping from the $k$th
symmetric power of $\Cal F^{(1)}$ onto $\Cal F^{(k)}$.
If $n+1 = \dim_{\Bbb R}(\Cal F^{(1)})$,
it follows that $\dim_{\Bbb R}(\Cal F^{(k)}) \le \binom{n+k}{k}$.

   (ii) Many examples of sequences of polynomial spaces have one more
property: the union $\bigcup_{k=0}^{\infty} \Cal F^{(k)}$
is dense in $L^2(\Omega,\sigma)$.

\bigskip

\proclaim{1.5.\ Definition} Let $(\Omega,\sigma)$ and 
$(\Cal F^{(k)})_{k \ge 0}$ be as in Definition~1.3.
A {\bf cubature formula of strength $k$ on $\Omega$}
is a pair $(X,W)$, consisting of a finite subset $X \subset \Omega$
and a weight $W : X \longrightarrow \Bbb R_+^*$, which is a
cubature formula for $\Cal F^{(k)}$.
In case $W(x) = \sigma(\Omega)/\abs{X}$ for all $x \in X$, 
the set $X$ is a {\bf geometrical
\footnote{
In particular cases, we replace \lq\lq geometrical\rq\rq \ by more
suggestive adjectives such as 
\lq\lq {\bf interval}\rq\rq \ in the situation 
of Example~1.10 with $\Omega$ an interval in $\Bbb R$, 
\lq\lq {\bf spherical}\rq\rq \ in the situation
of $\Omega = \Bbb S^{n-1}$, 
or \lq\lq {\bf Gaussian}\rq\rq \ if $\Omega$ is a real vector space 
and $\sigma$ a Gaussian measure. 
However, it should be kept in mind that the
\lq\lq Euclidean designs\rq\rq \ 
of \cite{NeuSe--88} and \cite{DelSe--89} (among others)
are of a different nature,
since the reference measure $\sigma$ 
on the Euclidean space $\Omega = \Bbb R^n$
depends in this case on the weighted set $(X,W)$
through some radial factor.
}
$k$-design} on $\Omega$. 
\endproclaim

\bigskip

\noindent {\bf 1.6.\ Organising questions.}
   (i) For standard examples of finite subsets $X$ of $\Omega$,
compute the largest strength $k$ 
for which $X$ is a geometrical $k$-design.

(ii) For $k$ given, 
find designs and cubature formulas $(X,W)$ of strength $k$
with $\abs{X}$ minimal. 
(More on this in Item~1.17, on tight spherical designs.)
In case $\Omega$ is moreover a metric space,
describe the distance set
$$
D_X \, = \, \left\{ c \in \mathopen]0,\infty\mathclose[ \ \big\vert  \
c = d(x,y) \quad \text{for some} \quad x,y \in X, x \ne y \right\} .
$$

   (iii) Dually, for given $N$, 
find designs and cubature formulas $(X,W)$ with
$\abs{X} = N$ of maximum strength.

   (iv)  Asymptotics and equidistribution.
Find sequences $\left( X^{(k)},W^{(k)} \right)_{k \ge 1}$,
where each $(X^{(k)},W^{(k)})$ is a cubature formula of strength $k$,
such that the sequence of measures $\left( \sigma^{(k)} \right)_{k \ge 1}$
converges to $\sigma$, where
$\sigma^{(k)} = \sum_{x \in X^{(k)}} W^{(k)}(x) \delta_x$
and where $\delta_x$ denotes the Dirac measure of
support $\{x\}$. 
Optimize in some sense the speed of convergence.

   When each $X^{(k)}$ is a geometrical $k$-design, this is the
standard question of $\sigma$-equidistri\-bution of a sequence of finite
subsets of a measure space; see \cite{Serre--97}.

\bigskip

  In the following proposition, $\ell^2(X,W)$ stands
for the Hilbert space of real-valued functions on $X$,
with scalar product defined by
$\langle \psi_1 \mid \psi_2 \rangle = \sum_{x \in X} W(x)\psi_1(x)\psi_2(x)$.

\proclaim{1.7.\ Proposition} 
Let $(\Omega,\sigma)$ be a finite measure space. 
Let $\Cal H$ be a finite dimensional space of functions 
in $L^2(\Omega,\sigma)$.
Assume that there exists a finite dimensional subspace $\Cal H '$
of $L^2(\Omega,\sigma)$  such that $\Cal H$ is linearly generated by
products $\varphi_1 \varphi_2$ with $\varphi_1,\varphi_2 \in \Cal H '$.
Let $X$ be a finite subset of $\Omega$
and let $W : X \longrightarrow \Bbb R_+^*$ be a weight.

Then $(X,W)$ is a cubature formula for $\Cal H$ if and only if
the restriction mapping
$$
\rho \, : \, \Cal H ' \longrightarrow \ell^2(X,W)
$$
is an isometry.
In particular, if $(X,W)$ is a cubature formula for $\Cal H$,
then
$$
\abs{X} \, \ge \, \dim_{\Bbb R}(\Cal H') .  \tag$*$
$$
\endproclaim

\demo{Proof} The condition for $\rho$ to be an isometry is
$$
\sum_{x \in X} W(x) \varphi_1(x)  \varphi_2(x)
\, = \, 
\int_{\Omega} \varphi_1(\omega)  \varphi_2(\omega)
d\sigma(\omega)
$$
for all $\varphi_1,\varphi_2 \in \Cal H '$.
By hypothesis on $\Cal H$ and $\Cal H'$, 
this condition is equivalent to
$$
\sum_{x \in X} W(x) \varphi(x) 
\, = \, 
\int_{\Omega} \varphi(\omega) d\sigma(\omega)
$$
for all $\varphi \in  \Cal H$.
$\square$
\enddemo

\bigskip

\proclaim{1.8.\ Definition}  Let $(\Omega,\sigma)$ and 
$(\Cal F^{(k)})_{k \ge 0}$ be as in Definition~1.3.
A cubature formula $(X,W)$ of strength $2l$
is {\bf tight} if equality holds in $(*)$ above,
with $\Cal H ' = \Cal F^{(l)}$ and $\Cal H = \Cal F^{(2l)}$.

   For a spherical cubature formula of strength $2l+1$,
\lq\lq tight\rq\rq \ is defined in~1.11.
\endproclaim

   Tight cubature formulas are rare: see Discussion~1.17. 

\bigskip

\noindent {\bf 1.9.\ Quadrature formulas on intervals.} 
Let $\Omega$ be a real interval
$(a,b)$, with $-\infty \le a < b \le \infty$
(the interval can be open or closed).
The spaces $\Cal F^{(k)}(a,b)$ of polynomial functions on $\Omega$ of
degree at most $k$ provide the canonical example for Definition~1.3.

Consider in particular the case $\Omega = [-1,1]$, with Lebesgue
measure, and an integer $l \ge 1$. 
Let first 
$$
X^{(l)} \, = \, \left\{ x_{1,l}, x_{2,l}, \hdots, x_{l,l} \right\}
\, \subset \, [-1,1]
$$
be an arbitrary subset of $l$ distinct points in $[-1,1]$.
Let $L_j(t) = \prod\frac{t-x_{i,l}}{x_{j,l} - x_{i,l}}$
(product over $i \in \{1,\hdots,l\}$, $i \ne j$)
be the corresponding {\it elementary Lagrange interpolation polynomials}.
Then $(X^{(l)},W^{(l)})$ is a \lq\lq cubature formula\rq\rq \
for $\Cal F^{(l-1)}(-1,1)$,
with the \lq\lq weight\rq\rq \ $W^{(l)}$ defined on $X^{(l)}$ by
$$
W^{(l)}\left(x_{j,l} \right) \, = \,
\int_{-1}^1 L_j(t) dt \qquad (1 \le j \le l) .
$$
But quotation marks are in order since the \lq\lq weight\rq\rq \ values
need not be positive; 
for example, in case $x_{j,l} = -1 + 2\frac{j-1}{l-1}$,
these {\bf Newton-Cotes formulas} have all weights positive if and only if
either $l \le 8$ or $l = 10$ (see, e.g., Chapter~2 of \cite{CroMi--84}).

Let now $X^{(l)}$ be precisely
the set of the roots of the Legendre polynomial of degree~$l$,
namely of
$P_l(t) = (-1)^l \frac{l!}{2l!}\frac{d^l}{dt^l}\left((1-t^2)^l\right)$.
The pair $(X^{(l)},W^{(l)})$ is a cubature formula for
$\Cal F^{(2l-1)}(-1,1)$,
known as a {\bf Gauss quadrature}
\footnote{
If $\Omega$ is an interval of $\Bbb R$, the term
\lq\lq quadrature\rq\rq \ is often used instead of \lq\lq cubature\rq\rq 
\ as in Definition~1.1.
}
(or sometimes {\it Gauss-Jacobi mecha\-nical quadrature});
this is the unique cubature formula  with $\le l$ points
for $\Cal F^{(2l-1)}(-1,1)$ and Lebesgue measure.
The weights $W^{(l)}(x_{j,l})$ of the Gauss formula
are strictly positive; indeed, since the polynomial $L_j$ defined above
have now their {\it squares} in $\Cal F^{(2l-1)}(-1,1)$, we have
$$
W^{(l)}(x_{j,l}) \, = \,
\sum_{k=1}^{l} W^{(l)}(x_{k,l}) L_j(x_{k,l})^2 \, = \, 
\int_{-1}^1 L_j(t)^2 dt \, > \, 0
$$
for $j \in \{1,\hdots,l\}$.
(See, e.g., \cite{Gauts--97}.)
\par

In general, there are analogous formulas
for other intervals of the real line and other measures 
with finite moments and with 
$\operatorname{dim}_{\Bbb R}\left(L^2((a,b),\sigma)\right) = \infty$.
This is a part of the theory of orthogonal polynomials;
see Section~3.4 in \cite{Szeg\"o--39}.
There are also related results in larger dimensions;
see Section 3.7 in \cite{DunXu--01}.  

\bigskip 

\noindent {\bf 1.10.\ Interval designs.} 
It is a particular case of 
a theorem of Seymour and Zaslavsky that 
interval designs (named \lq\lq averaging sets\rq\rq \ by these authors)
exist for any finite-dimensional space of continuous functions
on $\Omega = \mathopen]0,1\mathclose[$, with various measures;
see \cite{SeyZa--84} and \cite{Arias--88}.

An interval design for the space $\Cal F^{(2l-1)}(-1,1)$ 
of polynomial functions of degree at most~$2l-1$,
and Lebesgue measure,
requires strictly more than $\frac{1}{4}l^2$ points 
(this is a 1937 result of Bernstein for which we refer to
\cite{Krylo--62} and \cite{Korev--94}),  
and there exists a constant $c$ for which it is known that
$ck^{3}$ points suffice 
(Section 3.3 in \cite{Korev--94}, 
improving on Theorem 2 in \cite{Bajn--91a}).
For low values of $k$, as an answer to a question of Chebyshev (1874),
Bernstein has shown  that
there exist interval designs for $\Cal F^{(k)}(-1,1)$
with $\abs{X} = k$ 
if and only if $k \le 7$ or $k = 9$,
and there is an explicit construction for these values of $k$
(see Section~10.3 in \cite{Krylo--62}).

\bigskip

\noindent {\bf 1.11.\ Cubature formulas on spheres and tightness.} 
Let $\Omega = \Bbb S^{n-1}$ be the unit sphere 
in the Euclidean space $\Bbb R^n$, $n \ge 2$, and let
$\sigma$ denote the probability measure on $\Bbb S^{n-1}$ which is
invariant by the orthogonal group $O(n)$. 

   For spheres, the standard example for Definition~1.3 is given by
$$
\Bbb R = \Cal F^{(0)}(\Bbb S^{n-1}) \, \subset \, 
\Cal F^{(1)}(\Bbb S^{n-1}) \, \subset \, \cdots \, \subset \, 
\Cal F^{(k)}(\Bbb S^{n-1}) \, \subset \, \cdots
$$
where $\Cal F^{(k)}(\Bbb S^{n-1})$  denotes the space of restrictions to
the sphere of polynomial functions on $\Bbb R^n$ of degree at most $k$.
In this paper, when Definition~1.8 is particularized  
to spherical cubature formulas and designs, 
it is always with respect to this sequence of polynomial spaces;
thus, a cubature formula $(X,W)$ or strength $2l$ 
is {\bf tight} if
$$
\abs{X} \, = \, \dim_{\Bbb R}\big(\Cal F^{(l)}(\Bbb S^{n-1})\big)
\, = \, \binom{n+l-1}{n-1} + \binom{n+l-2}{n-1} .
$$

Consider the space $\Cal P^{(k)}(\Bbb S^{n-1})$ of
restrictions to $\Bbb S^{n-1}$ of polynomial functions on 
$\Bbb R^n$ which are homogeneous of degree $k$, and
the space $\Cal H^{(k)}(\Bbb S^{n-1})$ of
restrictions to $\Bbb S^{n-1}$ of polynomial functions on 
$\Bbb R^n$ which are homogeneous of degree $k$ and harmonic
\footnote{
A smooth function $\varphi$ on $\Bbb R^n$ is {\bf harmonic}
if $\sum_{j=1}^n \frac{\partial ^2 \varphi}{\partial x_j^2} = 0$.
}
\hskip-.1cm
. 
We have
$$
\aligned
\Cal F^{(k)}(\Bbb S^{n-1}) \, &= \, 
\Cal P^{(k)}(\Bbb S^{n-1}) \bigoplus \Cal P^{(k-1)}(\Bbb S^{n-1})
\\
\Cal P^{(k)}(\Bbb S^{n-1}) \, &= \,
\bigoplus_{j=0}^{[k/2]} \Cal H^{(k-2j)}(\Bbb S^{n-1}) 
\endaligned
$$
for all $k \ge 0$
(see e.g. Section IV.2 in \cite{SteWe--71}).
For reasons related to Proposition~3.2 and Remark~3.3 below, 
a cubature formula $(X,W)$ or strength $2l+1$ 
is {\bf tight} if
$$
\abs{X} \, = \, 2\dim_{\Bbb R}\big(\Cal P^{(l)}(\Bbb S^{n-1})\big)
\, = \, 2\binom{n+l-1}{n-1}.
$$
\par

Another example for Definition~1.3 is given by
$$
\Bbb R = \Cal P^{(0)}(\Bbb S^{n-1}) \, \subset \, 
\Cal P^{(2)}(\Bbb S^{n-1}) \, \subset \, \cdots \, \subset \, 
\Cal P^{(2l)}(\Bbb S^{n-1}) \, \subset \, \cdots \, 
$$
(observe that 
$\Cal P^{(k)}(\Bbb S^{n-1}) \subset \Cal P^{(k+2)}(\Bbb S^{n-1})$
for all $k \ge 0$,
since the restriction of $x_1^2 + \cdots + x_n^2$ to $\Bbb S^{n-1}$
is the constant $1$).
\par

A subset $X \subset \Bbb S^{n-1}$ is called {\bf antipodal}
if $-X = X$; such a $X$ can always be written (non-uniquely)
as a disjoint union of some $Y$ inside $X$ and of $-Y$.
A cubature formula $(X,W)$ on $\Bbb S^{n-1}$ is {\bf antipodal} 
if $X$ is antipodal and if $W(-x) = W(x)$ for all $x \in X$.
If $(X,W)$ is antipodal,
$\sum_{x \in X}W(x)\varphi(x) = 
\int_{\Bbb S^{n-1}}\varphi(\omega) d\omega = 0$
for any homogenous polynomial function $\varphi$ of odd degree;
it follows that $(X,W)$ is 
a cubature formula for $\Cal F^{(2l+1)}(\Bbb S^{n-1})$
if and only if it is 
a cubature formula for $\Cal P^{(2l)}(\Bbb S^{n-1})$.

\bigskip

\noindent {\bf 1.12.\ On strength values for spherical designs and
cubature formulas.} 
The following existence result goes back essentially to
the solution of Waring's problem by Hurwitz and Hilbert
(see \cite{Hurwi--08}, \cite{Hilbe--09},  page~722 in \cite{Dicks--19},
\cite{Ellis--71}, and \cite{Natha--96}):

   {\it for any $l \ge 0$, there exists a cubature formula
   $\left(Y^{(l)},W^{(l)}\right)$ for $\Cal P^{(2l)}(\Bbb S^{n-1})$
   with}
   $$
   \aligned
   \binom{n+l-1}{n-1} \, = \,  
   &\dim_{\Bbb R}\big(\Cal P^{(l)}(\Bbb S^{n-1})\big)  \, \le \, 
   \abs{ Y^{(l)} } \, \le \, 
   \\
   &\dim_{\Bbb R}\big(\Cal P^{(2l)}(\Bbb S^{n-1})\big) - 1 \, = \, 
   \binom{n+2l-1}{n-1}  - 1 .
   \endaligned
   $$
The lower bound is that of Proposition~1.7; 
the upper bound is a particular case of Theorem~2.8 below, 
and improves by $1$ the bound of Theorem~2.8 in \cite{LyuVa--93}.
The relation with Waring's problem is alluded to in Item~3.9;
for the application to Waring's problem, it is necessary
for the weights of the cubature formula to be {\it rational}
positive numbers; our proof of Theorem 2.8 does not show this,
and we refer to Chapter 3 of \cite{Natha--96}
for a complete proof.
 
It follows that 
{\it there exists an antipodal cubature formula
$\left(X^{(l)},W^{(l)}\right)$ of strength $2l+1$ 
with
$$
   2\binom{n+l-1}{n-1} \, \le \, 
   \abs{ X^{(l)} } \, \le \, 
   2\binom{n+2l-1}{n-1}  - 2 .
$$
}

However, the proof of the existence of 
$\left(Y^{(l)},W^{(l)}\right)$ 
does not provide nice constructions. 
For $n \ge 3$ and $k \ge 4$, we do not know many good
explicit cubature formulas of strength $k$ on $\Bbb S^{n-1}$.

   Note the relevant asymptotics. For $n$ fixed:
$$
\binom{n+l-1}{n-1} \, \approx \, \frac{l^{n-1}}{n!}
\qquad \text{if} \quad l \to \infty .
$$
For $l$ fixed:
$$
\binom{n+l-1}{n-1} \, \approx \, \frac{n^l}{l!}
\qquad \text{if} \quad n \to \infty .
$$
The lower bound
$\abs{X^{(l)}} \ge 2 \binom{n+l-1}{n-1}$
has been improved in some cases, 
in particular for $n$ fixed and $l \to \infty$ \cite{Yudin--97}.

\bigskip

   The general result of Seymour and Zaslavsky already quoted in Item~1.10
implies:

\smallskip

{\it for any $n \ge 2$ and $k \ge 0$, there exists a spherical
$k$-design $X \subset \Bbb S^{n-1}$.}

\smallskip

\par\noindent
The proof of \cite{SeyZa--84} is not constructive and
does not give any bound on the size
$\abs{X}$.

\bigskip

Let  $X$ be a non-empty finite subset of $\Bbb S^{n-1}$
and let $W$ be a weight on $X$
such that $\sum_{x \in X} W(x) = 1$. 

\smallskip

   (i) The pair $(X,W)$ is always a cubature formula of strength $0$.
It is a cubature formula of strength $1$ if and only if 
the weighted barycentre $\sum_{x \in X} W(x)x$ of $(X,W)$
coincides with the origin of $\Bbb R^n$.

(ii) If $(X,W)$ is a cubature formula of strength $2$,
then $X$ generates $\Bbb R^n$.
Indeed, if there exists a non-zero vector $\alpha$ orthogonal to $X$,
the function $\omega \longmapsto \langle \alpha \mid \omega \rangle^2$
has a non-zero integral on $\Bbb S^{n-1}$ but vanishes identically on $X$.

{\it Exercise.} Denote by $Z$ the {\it Gram matrix} of $X$,
defined by $Z_{x,y} = \langle x \mid y \rangle$ for $x,y \in X$,
and by $J$ the matrix with rows and columns indexed by $X$ 
and with all entries $1$.
Show that $X$ is a spherical $2$-design if and only if 
the three following conditions hold:
   $(ii_a)$ $Z_{x,x} = 1$ for all $x \in X$, 
   $(ii_b)$ $ZJ = 0$, 
   $(ii_c)$ $Z^2 = n^{-1}\abs{X}Z$. 
(Solution in Lemma~13.6.1 of \cite{Godsi--93}.)  

(iii) If $(X,W)$ is a cubature formula of strength $4$,
then $X$ cannot be a disjoint union of two orthogonal sets; 
more generally, 
in case there exist two vectors $\alpha,\beta \in \Bbb S^{n-1}$
such that $X \subset \alpha^{\perp} \cup \beta^{\perp}$,
the function 
$\omega \longmapsto
\langle \alpha \mid \omega \rangle^2
\langle \beta \mid \omega \rangle^2$
has a non-zero integral on $\Bbb S^{n-1}$ but vanishes identically on $X$.

(iv) If $(X,W)$ is a cubature formula on $\Bbb S^{n-1}$ of strength $2l$
for some $l \ge 1$, then the set
$P^{(l)}(X) = \left\{ 
\omega \longrightarrow \langle \omega \mid x \rangle^l
\right\}_{x \in X}$
linearly generates $\Cal P^{(l)}(\Bbb S^{n-1})$.
% Ref de Thierry Vust :
% Lemma 1.5.1 dans 
% D.J. Benson, Polynomial invariants of finite groups
% Cambridge Univ. Press 1993, GE 20/296
Indeed, assume that some $\varphi \in \Cal P^{(l)}(\Bbb S^{n-1})$
is orthogonal to $P^{(l)}(X)$.
If $\hat \varphi \in \Cal P^{(l)}(\Bbb S^{n-1})$ is defined by
$$
\hat \varphi (u) \, = \, 
\int_{\Bbb S^{n-1}} \langle \omega \mid u \rangle^l
\varphi(\omega) d \sigma (\omega)
\quad \forall \ u \in \Bbb R^n ,
$$
the hypothesis means that the restriction of $\hat \varphi$
to $X$ vanishes;  this implies that $\hat \varphi = 0$
by Proposition~1.7. 
But $\hat \varphi = 0$ implies $\varphi = 0$ since 
$\left\{  \omega \longrightarrow \langle \omega 
\mid u \rangle^l \right\}_{u \in \Bbb S^{n-1}}$
linearly generates $\Cal P^{(l)}(\Bbb S^{n-1})$.
With a terminology borrowed from the theory of lattices (see
Reminder~1.15), the special case of Property (iv) for $l = 2$ 
states that cubature formulas of strength $4$ on $\Bbb S^{n-1}$ 
provide perfect sets.

\bigskip 

   For cubature formulas on spheres, 
the following equivalences are useful. 
Since many natural examples are provided by 
intersections of lattices with spheres of varius radii (see~1.16), 
we find it useful to consider for each $\rho > 0$ 
the sphere $\rho\Bbb S^{n-1}$ of radius $\rho$ and centre the
origin in $\Bbb R^n$; 
we denote again by $\sigma$ the $O(n)$-invariant
probability measure on $\rho\Bbb S^{n-1}$.

\bigskip

\proclaim{1.13.\ Proposition} 
Consider integers $n \ge 2$ and $k \ge 0$, 
a positive number $\rho$, 
a non-empty finite subset $X$ of $\rho\Bbb S^{n-1}$,
and a weight $W : X \longrightarrow \Bbb R_+^*$
such that $\sum_{x \in X}W(x) = 1$. 
The following conditions are equivalent
\footnote{
Even if there are canonical identifications between
$\Cal P^{(k)}(\Bbb R^n)$, $\Cal P^{(k)}(\Bbb S^{n-1})$, 
and $\Cal P^{(k)}(\rho\Bbb S^{n-1})$,
see Example~1.11, we choose here the first notation. 
}
\hskip-.1cm
.
\roster

\item"(i)" 
$
\sum_{x \in X} W(x) \varphi (x) =
\int_{\Bbb \rho S^{n-1}} \varphi(\omega) \, d\sigma (\omega) 
$
for all $\varphi \in \Cal P^{(j)}(\Bbb R^n)$ 
and $j \in \{0,1,\hdots,k\}$, 
namely $(X,W)$ is a cubature formula of strength $k$
on $\rho\Bbb S^{n-1}$. 
\par

\item"(ii)" 
$\sum_{x \in X} W(x) \varphi(x) = 
\sum_{x \in X} W(x) \varphi(g x)$
for all $\varphi \in \Cal P^{(j)}(\Bbb R^n)$, 
$j \in \{0,1,\hdots,k\}$,
and $g \in O(n)$.
\par

\item"(iii)" 
$\sum_{x \in X} W(x) \varphi(x) = 0$ 
for all $\varphi \in \Cal H^{(j)}(\Bbb R^n)$
and $j \in \{1,\hdots,k\}$. 
\par

\item"(iv)" 
If $l$ is the largest integer such that $2l \le k$, 
there exists a constant $c_{2l}$ such that 
$$
(iv)_a \hskip1truecm
\sum_{x \in X} W(x) \langle x \mid u \rangle ^{2l} \, = \, 
 c_{2l} \, \rho^{2l} \langle u \mid u \rangle ^{l} 
\qquad \forall \ u \in \Bbb R^n , 
$$
and, 
if $l'$ is the largest integer such that $2l'+1 \le k$,
then 
$$
(iv)_b \hskip2truecm
\sum_{x \in X} W(x) \langle x \mid u \rangle ^{2l'+1} \, = \, 0 
\qquad \forall \ u \in \Bbb R^n . \hskip1truecm
$$

\item"(v)"
Condition $(iv)_a$ holds for any even integer $2l \le k$
and Condition $(iv)_b$ holds for any odd integer $2l'+1 \le k$.

\endroster
Moreover, the constant in $(iv)_a$ is
$$
   c_{2l} \, = \, 
       \frac{1 \times  3 \times  5 \times \cdots \times (2l-1)}
      {n \, (n+2) \, (n+4) \, \cdots \, (n+2l-2)} 
$$
if $l > 0$, and $c_0 = 1$ if $l = 0$.
\endproclaim

\demo{Proof} The proof of Theorem~3.2 in \cite{VenMa--01}, 
for spherical designs, readily carries over to cubature formulas.
$\square$
\enddemo

\demo{Remarks} If $(X,W)$ is antipodal, then Equality $(iv)_b$ is
automatically satisfied.

Besides the conditions of Proposition 1.13, 
there are many other equivalent ones; 
see \cite{DeGoS--77} and \cite{GoeSe--79}, as well as Item~2.14.
\enddemo

\bigskip

\noindent {\bf 1.14.\ Examples of spherical designs:
orbits of finite groups, inductive cons\-tructions,
distance-regular graphs, and contact points of John's ellipsoids.}

\medskip\noindent
{\it Group orbits.}
\smallskip

   It is easy to show that any orbit in $\Bbb S^{n-1}$ 
of any irreducible finite subgroup of $O(n)$ 
is a $2$-spherical design; see \cite{Banna--79},
but the result can already be found in a 1940 paper by
R. Brauer and H.S.M. Coxeter
(Theorem 3.6.6 in \cite{Marti--03}).
Moreover, a $(2l)$-spherical design which is antipodal  
is also a $(2l+1)$-spherical design.  
In parti\-cular, 
if $(e_1,\hdots,e_n)$ is the canonical basis of $\Bbb R^n$, 
then the set $\{\pm e_1, \hdots, \pm e_n\}$ is a 
spherical $3$-design on $\Bbb S^{n-1}$ of cardinality~$2n$
(hence a tight $3$-design); 
the set $\{\pm e_1 \pm e_2 \pm \cdots \pm e_n\}$ 
(all choices of signs)
is another spherical $3$-design, of cardinality~$2^n$, 
and therefore non tight when $n \ge 3$; 
neither of these is a spherical $4$-design. 

   Orbits of groups generated by reflections provide interesting
spherical designs. In particular irreducible root systems of type
$A_n, D_n, E_n$ are spherical $3$-designs
(see Chapter~V, \S~6, Number~2 in \cite{Bourb--68}).
Moreover, roots systems of type $A_2, D_4, E_6, E_7$ provide spherical
$5$-designs and the root system of type $E_8$ provides a spherical
$7$-design. 

   For $n \ge 3$, there is a maximal strength 
$k_{\operatorname{max}}(n)$
for spherical designs in $\Bbb S^{n-1}$ which are orbits
of finite subgroups of $O(n)$.
It is moreover conjectured that 
$\sup_{n \ge 3}k_{\operatorname{max}}(n) < \infty$;
the maximal value of $k$ we know is $k = 19$ and
occurs in dimension $n = 4$ (see~3.8 below).
For more on spherical designs which are orbits of finite groups,
see \cite{Banna--79} and \cite{HarPa--04}.

\medskip\noindent
{\it Inductive constructions.}
\smallskip

   Spherical designs can be constructed inductively as follows.
On $\Bbb S^1$, the vertices of a regular $N$-gon constitute a
$k$-design if and only if $N \ge k+1$.
For $n \ge 3$ and some given $k$, assume that we have a
spherical $k$-design $Y$ on $\Bbb S^{n-2}$,
as well as an interval design $Z \subset \mathopen]-1,1\mathclose[$
for the space $\Cal F^{(k)}(-1,1)$ of polynomial functions 
of degrees at most $k$ 
and for the measure $d\sigma(t) = (1-t^2)^{(n-3)/2}dt$.
Then 
$$
X \, = \, \left\{ \left( \sqrt{1-\norm{z}^2}y \, , \, z \right) \in \Bbb
S^{n-1} 
\ \Big\vert \
y \in Y \quad
\text{and} \quad z \in Z \right\}
$$
is a spherical $k$-design in $\Bbb S^{n-1}$.
Denote by $M'_n(k)$ the smallest integer such that,
for any $N \ge M'_n(k)$, there exists a spherical $k$-design
of size $N$ on $\Bbb S^{n-1}$. With the construction above,
it can be shown that
$$
M'_n(k) \, \le \, O\left( k^{\frac{n(n+1)}{2}-2} \right) .
$$
It is conjectured that $\frac{n(n+1)}{2}-1$ can be reduced to
$\frac{n(n-1)}{2}$.
For all this, see \cite{Bajn--91b} and other papers by Bajnok.

   In case $n = 3$, it is a result of \cite{KorMe--93} that,
for all $k \ge 0$, there exists a spherical $k$-design of size
$O(k^3)$. 

   Constructions of spherical designs appear also in \cite{Mimur--90},
\cite{Bajno--92}, and \cite{Bajno--98}.

\medskip\noindent
{\it Distance regular graphs.}
\smallskip

Consider a distance-regular graph $\Gamma$,
with vertex set $V(\Gamma)$, of valency $k$.
Let $\theta$ be an eigenvalue of $\Gamma$, $\theta \ne k$,
of some multiplicity $m$;
let $p$ denote the orthogonal projection of
the space of functions $V(\Gamma) \longrightarrow \Bbb R$
to the $\theta$-eigenspace of the adjacency matrix of $\Gamma$.
Then $p(V(\Gamma))$ is a spherical $2$-design in
$\rho \Bbb S^{m-1}$ for the appropriate radius $\rho$
(Corollary~13.6.2 in \cite{Godsi--93}).

\medskip\noindent
{\it Contact points of John's ellipsoids.}
\smallskip

   Let $K$ be a convex body in $n$-space  
such that the maximal volume ellipsoid in $K$ is the unit n-ball $B$. 
It is a theorem of John that there exist
a finite subset $X$ of $K \cap \partial B$ and
a positive weight function $W$ on $X$
such that $(X,W)$ is a cubature formula 
for $\partial B = \Bbb S^{n-1}$ of strength~$2$,
indeed of strength $3$ if $K$ is assumed to be symmetric;
see Lecture 3 in \cite{Ball--97}.
It could be interesting to investigate systematically
polyhedra $K$ such that $K \cap \partial B$ is a spherical $k$-design,
or such that there exists a weight $W$ for which 
$(K \cap \partial B,W)$ is a cubature formula of strength $k$,
for various values of $k$.

\bigskip

\noindent {\bf 1.15.\ A reminder on lattices.}
Let $V$ be a Euclidean space of dimension $n \ge 1$.
A {\bf lattice} is a discrete subgroup $\Lambda$ of $V$
which generates $V$ as a vector space,
and $n$ is the {\bf rank} of $\Lambda$.
If $\Lambda \subset V$ is a lattice, so is its
{\bf dual} $\Lambda^* = \left\{ x \in V \mid 
\langle x \mid \Lambda \rangle \subset \Bbb Z \right\}$.
A lattice $\Lambda$ is {\bf integral}
if $\langle \Lambda \mid \Lambda \rangle \subset \Bbb Z$,
equivalently if $\Lambda \subset \Lambda^*$.
A lattice $\Lambda$ is {\bf unimodular}
if $\operatorname{Vol}(V/\Lambda) = 1$;
an integral lattice is unimodular
if and only if $\Lambda = \Lambda^*$.
A lattice $\Lambda$ is {\bf even}
if $\langle \lambda \mid \lambda \rangle \in 2\Bbb Z$
for every $\lambda \in \Lambda$;
an {\bf odd} lattice is an integral lattice which is not even.
Two lattices $\Lambda \subset V$, $\Lambda' \subset V'$
are {\bf equivalent} if there exists 
an isometry $g$ from $V$ onto $V'$ such that
$g(\Lambda) = \Lambda'$.
There is a natural notion of orthogonal direct sum
of lattices, and a lattice is {\bf irreducible}
if it is not equivalent to a direct sum
$\Lambda \oplus \Lambda '$ with $\Lambda,\Lambda ' \ne 0$.

Let $\Lambda \subset V$ be an integral lattice;
for $m \ge 1$, we denote by 
$$
\Lambda_m \, = \, \left\{ \lambda \in \Lambda \mid
\langle \lambda \mid \lambda \rangle = m \right\}
$$
the {\bf shell} of radius $m$
(namely of radius $\sqrt m$ in the usual sense
of Euclidean geometry).
For any $\lambda \in \Lambda_1$, there is 
an integral lattice $\Lambda'$ in the orthogonal subspace
$\lambda^{\perp}$ such that 
$\Lambda = \Bbb Z \lambda \oplus \Lambda '$;
in particular, $\Lambda_1 = \emptyset$
for an irreducible lattice of rank at least $2$.

A standard example of an odd unimodular lattice is the
{\bf cubical lattice} $\Bbb Z^n \subset \Bbb R^n$.
If $n \le 11$, any odd unimodular lattice is 
equivalent to one of these.

   Let $n$ be a multiple of $4$. 
In the standard Euclidean space $\Bbb R^n$,
define the lattice 
$D_n = \left\{ \lambda \in \Bbb Z^n \mid 
\sum_{i=1}^n \lambda_i \equiv 0 \pmod{2} \right\}$
and the {\bf Witt lattice}
$\Gamma_n = D_n \cup 
\left( (\frac{1}{2}, \hdots, \frac{1}{2}) + D_n\right)$.
Then $\Gamma_n$ is integral and unimodular;
moreover, it is even if and only if $n$ is a multiple of $8$.
Even unimodular lattices exist only 
in dimensions $n \equiv 0 \pmod{8}$.
If $n = 8$, any even unimodular lattice is equivalent to $\Gamma_8$,
also known as the {\bf root lattice of type $E_8$}
or the {\bf Korkine-Zolotareff lattice}.
If $n = 16$, any even unimodular lattice is equivalent to
either $\Gamma_{16}$ (which is irreducible)
or to $\Gamma_8 \oplus \Gamma_8$.
If $n = 24$, there are $24$ equivalence classes of even unimodular lattices
(Niemeier's classification, 1968);
the most famous of them is the lattice discovered by {\bf Leech} (1964),
which is the only even unimodular lattice $\Lambda$ in dimension $n \le 31$
such that $\Lambda_2 = \emptyset$,
and which has a remarkably high density 
(see, e.g., Table~1.5 in Chapter 1 of \cite{ConSl--99}).
If $n = 32$, the number of equivalence classes of even unimodular lattices
is larger than $8 \times 10^7$.

There is a classification of integral unimodular lattices of small rank.
For $n \le 16$, Kneser (1957) has shown
that the only irreducible integral unimodular lattices
are $\Bbb Z$, $\Bbb \Gamma_n$ for $n = 8, 12, 16$,
and three odd lattices in dimensions $14$, $15$, and $16$.
In particular, if $9 \le n \le 11$, 
any odd unimodular lattice is equivalent to 
either $\Bbb Z^n$ or $\Gamma_8 \oplus \Bbb Z^{n-8}$.
Integral unimodular lattices have been later classified for
$n \le 23$ (Conway and Sloane, 1982)
and $n \le 25$ (Borcherds, 1984).
For $n \le 24$, see Chapters 16 and 17 of \cite{ConSl--99};
minor corrections to previous tables are given in \cite{Bache--97}.

  Let $\Lambda \subset V$ be a lattice 
(possibly neither integral nor unimodular).
For $r > 0$, let $B_r$ denote the ball or radius $r$ 
centred at the origin of $V$.
The {\bf sphere packing}  associated to $\Lambda$ is
$\bigcup_{\lambda \in \Lambda}(\lambda + B_r)$,
where $r$ denotes the largest real number
such that the balls $\lambda + B_r$
have disjoint interiors ($\lambda \in \Lambda$);
the {\bf density} of $\Lambda$ is the number
$$
\lim_{R \to \infty}
\operatorname{Vol} \Big( 
    \big(\bigcup_{\lambda \in \Lambda}(\lambda + B_r)\big) \cap B_R \Big)
\Big/ 
\operatorname{Vol}(B_R) 
\, = \,
\frac{\operatorname{Vol}(B_r)}{\operatorname{Vol}(V/\Lambda)} .
$$
The lattice $\Lambda$ is {\bf extreme} 
if this density is a local maximum
in the space of all lattices of dimension $n$;
since density is homothety-invariant, 
a unimodular lattice is extreme if its density
is a local maximum in the space of  
all unimodular lattices of dimension $n$, 
a space which can be identified with the double coset space
$SL(n,\Bbb Z) \setminus SL(n,\Bbb R) / O(n)$.
   Denote by 
$$
\Lambda_{\operatorname{short}} \, = \, 
\left\{ 
\mu \in \Lambda  \ , \ \mu \ne 0  \mid
\langle \mu \mid \mu \rangle \le
\langle \lambda \mid \lambda \rangle 
\quad \text{for all} \quad 
\lambda \in \Lambda \ , \ \lambda \ne 0 
\right\}
$$
the shell of {\bf short vectors} in $\Lambda$;
the lattice  $\Lambda$ is {\bf eutactic} if there exists a weight
$W~: \Lambda_{\operatorname{short}} \longrightarrow \Bbb R_+^*$
such that $(\Lambda_{\operatorname{short}},W)$ is a spherical cubature
formula of strength $3$,
and {\bf perfect} if the set
$\left\{ \omega \longmapsto \langle \omega \mid \lambda \rangle^2
\mid
\lambda \in \Lambda_{\operatorname{short}} \right\}$
linearly generates the space of
homogeneous polynomial of degree $2$ on $V$.
The first result involving both lattices and cubature formulas is the
following theorem of {\bf Voronoi} (1908):
\par\centerline{
{\it a lattice is extreme if and only if it is both eutactic and perfect.}
}
\par\noindent
For $\Lambda$ to be extreme, it is sufficient that
$\Lambda$ is {\bf strongly perfect}, which means that
$\Lambda_{\operatorname{short}}$ is a spherical $5$-design
(Theorem 6.4 in \cite{VenMa--01}).
Strongly perfect lattices have been classified in dimensions $n \le 11$,
where they occur in dimensions $n =1, 2, 4, 6, 7, 8, 10$ only
\cite{VenMa--01}, \cite{NebVe--00};
other examples occur in \cite{BacVe--01}.

   Standard references on lattices include 
\cite{ConSl--99}, 
\cite{Ebeli--94},
\cite{Marti--03}, 
\cite{MilHu--73},
\cite{Serre--70},
and \cite{VenMa--01}.

\bigskip

\noindent {\bf 1.16.\ Examples of spherical designs:
lattice designs.} 
   Let $\Lambda$ be a lattice in $\Bbb R^n$. 
Any non-empty shell $\Lambda_m$
provides a spherical design of some strength
in the sphere $\sqrt{m}\Bbb S^{n-1}$.
This connection goes back to \cite{Venko--84}. 
See \cite{Pache} for many examples; 
let us indicate here a sample.

   If $\Lambda \subset \Bbb R^8$ is a root lattice of type $E_8$,
then $\Lambda_{2m}$ is  a spherical $7$-design for any $m \ge 1$. 
Moreover, it can be shown that
a shell $\Lambda_{2m}$ is a $8$-design if and only if
the Ramanujan coefficient $\tau (m)$ of the modular form
$$
\aligned
\Delta_{24}(z) \, &= \, q^2 \prod_{m=1}^{\infty}(1-q^{2m})^{24}
\, = \, \sum_{m=1}^{\infty} \tau(m) q^{2m}  \\
\, &= \, q^2 - 24 q^4 + 256 q^6 - 1472 q^8 + \cdots 
\qquad (q = e^{i\pi z})
\endaligned
$$
is zero. (See Example~3.10 below for a second appearance of 
this modular form $\Delta_{24}$ of weight $12$.)
Now it is a famous conjecture of Lehmer \cite{Lehme--47} 
that $\tau(m) \ne 0$ for all $m \ge 1$.
If $\Gamma \subset \Bbb R^{16}$ is a Witt lattice
(an irreducible even unimodular lattice, 
uniquely defined up to isometry by these properties
in dimension $16$),
then $\Gamma_{2m}$ is a spherical $3$-design for any $m \ge 1$,
and the condition for one of these shells to be a $4$-design
happens to be again the vanishing 
of the corresponding Ramanujan coefficients.
The same holds for the reducible even unimodular lattice 
$\Lambda \oplus \Lambda \in \Bbb R^{16}$. 
Consequently, the following four claims are simultaneously 
true or not true:
{\it 
\roster
\item"(i)" Lehmer's conjecture holds, namely $\tau (m) \ne 0$ 
   for all $m \ge 1$;
\item"(ii)" no shell of the root lattice $\Lambda$ of type $E_8$ is a 
   spherical $8$-design;
\item"(iii)" no shell of the lattice 
   $\Lambda \oplus \Lambda \in \Bbb R^{16}$
   is a spherical $4$-design;
\item"(iv)" no shell of the Witt lattice 
   $\Gamma \in \Bbb R^{16}$
   is a spherical $4$-design.
\endroster
}
Claim (i) has been checked for $m \le 10^{15}$
\cite{Serre--85}.

  It is easy to formulate other equivalences of the same kind. For
example, any shell of the Leech lattice~$L$ (see 1.15) 
is a spherical $11$-design, and the condition for $L_{2m}$ ($m \ge 2$)
to be a $12$-design is that the $m$th coefficient
$\mu(m) = \sum_{j=1}^{m-1}\tau(j)\tau(m-j)$ of $(\Delta_{24})^2$
vanishes. Thus

\centerline{{\it
$\mu_m \ne 0$ for any $m \ge 2$ if and only if $L_{2m}$ is not a
spherical $12$-design 
}}

\noindent
and, moreover, we have checked that $\mu_m \ne 0$ for $m \le 1200$.
We do not know of any other lattice with a shell 
which is a spherical design of strength $k \ge 11$.

   Similarly, let now $\Lambda$ denote a Niemeier lattice, 
namely an even unimodular lattice in dimension $24$ 
with $\Lambda_2 \ne \emptyset$
(up to isometry, there are  $23$ such lattices). 
Denote by $Q$ the modular form of
weight $4$ defined by
$$
\aligned
Q(q) \, &= \, 
\left( \theta_3(q) \right)^8 - 
      \frac{1}{16} \left( \theta_2(q)\theta_4(q) \right)^4 \\
&= \,
\bigg( \sum_{m \in \Bbb Z} q^{m^2} \bigg)^8
-
\bigg( \sum_{m \in 1/2 + \Bbb Z} q^{m^2} \bigg)^4
\bigg( \sum_{m \in \Bbb Z} (-q)^{m^2} \bigg)^4 \\
&= \, 1 + 240q^2 + 2160 q^4 + \cdots 
\endaligned
$$
and set
$$
Q(q)\Delta_{24}(q) \, = \, \sum_{m=1}^{\infty} \nu_m q^m .
$$
Then
\par
\centerline{{\it
$\nu_m \ne 0$  if and only if $\Lambda_{2m}$ is not a
spherical $4$-design ($m \ge 1$)
}}

\noindent
and, moreover, we have checked that $\nu_m \ne 0$ for $m \le 1200$.
(Observe that the condition is the same 
for any of the $23$ Niemeier lattices.)

   For one more example in this class, consider the cubical lattice
$\Bbb Z^n$. The criterion stated in 1.14 in terms of irreducible
representations of finite subgroups of $O(n)$
shows that all non-empty shells are spherical $3$-designs.
Moreover, it can be shown that
there are two classes of \lq\lq special shells\rq\rq \ 
which are spherical $5$-designs, namely
$(\Bbb Z^4)_m$ for $m = 2a$
and $(\Bbb Z^7)_m$ for $m = 4^b(8a+3)$.
Let us restrict now for brevity to $n \ge 8$,
and denote by $\Theta^{[n]}$ the modular form of weight $4 + n/2$   
defined by 
$$
\Theta^{[n]}(q) \, = \,  
\frac{1}{16} \left( \theta_2(q)\theta_4(q) \right)^4 
\left( \theta_3(q) \right)^n
\, = \, \sum_{m=1}^{\infty} \kappa^{[n]}_m q^n .
$$
Then

\centerline{{\it
$\kappa^{[n]}_m \ne 0$ 
if and only if $\left(\Bbb Z^n\right)_{m}$
is not a spherical $4$-design ($m \ge 1$) 
}}

\noindent
and, moreover, we have checked that these hold 
for all $n \ge 8$ and $m \le 1200$.   
(See \cite{Pache}, which contains a discussion
including $1 \le n \le 7$.)

\bigskip

\noindent {\bf 1.17.\ Tight spherical designs 
on spheres of dimension $n-1 \ge 2$.}
(See \cite{BaMuV--04} and references there.)

\medskip\noindent
{\it Even strengths.}
\smallskip

Tight spherical $(2l)$-designs do not exist when $2l \ge 6$. 

Tight spherical $4$-designs in $\Bbb S^{n-1}$ are of size
$\operatorname{dim}_{\Bbb R}(\Cal F^{(2)}(\Bbb S^{n-1}))
= \frac{1}{2}n(n+3)$.
They cannot exist unless  $n$ is of the form $(2m+1)^2 - 3$. 
If $m = 1$ or $m = 2$ examples are known, 
and known to be unique up to isometry;
they are respectively of size $27$ in $\Bbb S^5$  
and size $275$ in $\Bbb S^{21}$ (see below).
If $m = 3$ and $m = 4$ 
(and infinitely many larger values), 
non-existence has been proved. 

Tight spherical $2$-designs exist in all dimensions, and are regular
simplices. 

\medskip\noindent
{\it Odd strengths.}
\smallskip

For $l \ge 0$, a tight spherical $(2l+1)$-design is necessarily
antipodal, by Theorem~5.12 in \cite{DeGoS--77}.

Tight spherical $(2l+1)$-designs do not exist when $2l+1 \ge 9$,
up to one exception  (which is unique up to isometry):  
the $196 \, 560$ short vectors 
of a Leech lattice which provide 
(after dividing all vectors by $2$)
an $11$-design in $\Bbb S^{23}$. 
Observe that   
$$
196560 \, = \, 2\binom{28}{5} \, = \,
2\operatorname{dim}_{\Bbb R}(\Cal P^{(5)}(\Bbb S^{23})) .
$$

Tight spherical $7$-designs in $\Bbb S^{n-1}$
are of size $\frac{1}{3}n(n+1)(n+2)$.
They cannot exist unless
$n$ is of the form $3m^2-4$.
If $m = 2$ or $m = 3$, examples are known,
and known to be unique up to isometry.
(Up to homothety, 
they are respectively a root system of type $E_8$ and 
the short vertices of the unimodular integral lattice denoted by
$O_{23}$ in \cite{VenMa--01}
\footnote{
Let $L$ be a Leech lattice. Let first $e \in L$ be a short vector,
with $\langle e \mid e \rangle = 4$.
Denote by $p$ the orthogonal projection of $\Bbb R^{24}$
onto $e^{\perp}$ and set
$L'_e = \{x \in L \mid \langle e \mid x \rangle \equiv 0
\pmod{2} \}$.
Then $O_{23} = p(L'_e)$ is a unimodular integral lattice with 
$\min\{\langle x \mid x \rangle \mid x \in O_{23} , \quad x \ne 0 \}
= 3$
and $\abs{ \{x \in O_{23} \mid \langle x \mid x \rangle = 3 \} }
= 4600$. 
Let then $f \in L$ be a vector such that
$\langle f \mid f \rangle = 6$. 
Then $M_{23} = L \cap f^{\perp}$ is an integral lattice with 
$\min\{\langle x \mid x \rangle \mid x \in M_{23} , \  x \ne 0 \} = 4$;
if $M_{23}^*$ denote its dual lattice, then $[M_{23}^* : M_{23}] = 4$,
and there are $552$ short vectors in $M_{23}^*$.
}
\hskip-.1cm
.)
If $m = 4$ and $m = 5$ 
(and infinitely many larger values),
non-existence has been proved.

Tight spherical $5$-designs in $\Bbb S^{n-1}$
are of size $n(n+1)$.
They cannot exist unless 
$n$ is either $3$ or of the form $(2m+1)^2 - 2$. 
If $n = 3$, or $m = 1$, or $m = 2$, examples are known, 
and known to be unique up to isometry. 
(They are respectively a regular icosahedron and,
up to homothety, the short vectors of
a lattice which is dual to a root lattice of type $E_7$ and 
the short vectors of a lattice 
of type $M_{23}^*$, with the notation of \cite{VenMa--01}.) 
If $m = 3$ and $m = 4$ 
(and infinitely many larger values),
non-existence has been proved. 

Let $X \subset \Bbb S^{n-1}$ be a tight spherical $5$-design.
It is a consequence of Proposition~3.4 that 
$\langle x \mid y \rangle = \pm \alpha$ 
for any $x,y \in X$ with $x \ne \pm y$,
where $\alpha = 1/\sqrt{n+2}$.
Choose $e \in X$ and set 
$X_0 = \{x \in X \mid \langle x \mid e \rangle = \alpha \}$.
Then $X_0$ is a tight spherical $4$-design in 
the sphere $\sqrt{1-\alpha^2}\Bbb S^{n-2}$
centred at $\alpha e$ in the affine hyperplane
$\left\{ x \in \Bbb R^n \mid \langle x \mid e \rangle = \alpha \right\}$,
by Theorem~8.2 in \cite{DeGoS--77}.

Tight spherical $3$-designs exist in all dimensions, 
and are of the form $\{\pm e_1, \hdots, \pm e_n\}$,
where $\{e_1, \hdots, e_n\}$ is an orthonormal basis of $\Bbb R^n$.
Tight spherical $1$-designs are of the form $\{\pm e_1\}$.

\medskip

For results on tight designs on other spaces $(\Omega,\sigma)$,
see \cite{Baba1} and \cite{Baba2}.

\bigskip

\noindent {\bf 1.18.\ Designs with few points on $\Bbb S^2$.} 

   Given an integer $N \ge 1$, it is a natural question to ask what is the
largest integer $k_N \ge 0$ for which there exists a spherical
$k_N$-design with $N$ points on a sphere, say here on $\Bbb S^2$.
We report now some answers to this question, from \cite{HarSl--96}.
$$
\matrix
N = 1 & k_1 = 0 &   && \text{single point} \\
N = 2 & k_2 = 1 &   && \text{pair of antipodal points} \\
N = 3 & k_3 = 1 &   && \text{equatorial equilateral triangle} \\
N = 4 & k_4 = 2 &   && \text{regular tetrahedron (tight)} \\
N = 5 & k_5 = 1 &   && \text{} \\
N = 6 & k_6 = 3 &   && \text{regular octahedron (tight)} \\
N = 7 & k_7 = 2 &   && \\
N = 8 & k_8 = 3 &  && \text{cube (non-tight)} \\
N = 9 & k_9 = 2 &   && \\
N = 10 & k_{10} = 3 &   && \\
N = 11 & k_{11} = 3 &   && \\
N = 12 & k_{12} = 5 &   && \text{regular icosahedron (tight)} \\
13 \le N \le 23 & k_N \in \{3,4,5\} &   &&  \\
N = 20 & k_{20} = 5 & && \text{regular dodecahedron (non tight)} \\
N = 24 & k_{24} = 7 & && \text{improved snub cube} \\
\cdots & \cdots & && \\
N = 60 & k_{60} = 10 &   &&
\endmatrix
$$

  {\it Remarks.} (i) The table shows that $k_N$ is {\it not} monotonic 
as a function of $N$.

   (ii) We know from Propositions~1.7 and~3.2 that, on $\Bbb S^2$,
any spherical $(2l)$-design has size $N \ge (l+1)^2$
and any spherical $(2l+1)$-design has size $N \ge (l+2)(l+1)$.
The table above shows that these bounds are not sharp
unless $2l \in \{0,2\}$ or $2l+1 \in \{1,3,5\}$.
For example, a $4$-design has minimal size $12 (> 9)$ 
and a $7$-design has minimal size $24 (> 20)$.

   (iii) The table refers to isolated examples; here is a continuous
family, from \cite{HarSl--92}.  
Distribute the $12$ points of a regular icosahedron 
into two  poles $N$ and $S$, 
and two sets $P, Q$ of $5$ points each in planes orthogonal 
to the diameter joining $N$ to $S$. Let
$X_{\theta}$ denote the union of the poles, the set $P$, and the image of
the set $Q$ by a rotation of angle $\theta$ fixing $N$ and $S$. If $\theta$
is not a multiple of $2\pi / 5$, then $X_{\theta}$ is a spherical
$4$-design which is not a $5$-design. 

   (iv) The Archimedes' regular snub cube
\footnote{
See for example
http://mathworld.wolfram.com/SnubCube.html
} 
($24$ vertices) is a spherical $3$-design (not a $4$-design), 
while the improved snub cube 
reported to in the table is indeed  a spherical $7$-design.

   (v) The regular truncated icosahedron (= soccer ball, $60$ vertices)
is a spherical $5$-design (not a $6$-design). 
The \lq\lq improved soccer ball\rq\rq \ of \cite{GoSe--81a/b}, 
which is almost indistinguishable from the regular one
with the naked eye, is a $9$-design. 
The spherical $10$-design of size $60$ which appears in
\cite{HarSl--96} is quite different.

   (vi) The same paper shows a $9$-design with $N = 48$ points; 
indeed $k_{48} = 9$ and $k_N < 9$ 
whenever $N < 48$ or $N \in \{49, 51, 53\}$.
On the other hand, there exists a cubature formula on $\Bbb S^2$
of strength $9$ and size $32$ 
(Section 5 in \cite{GoSe--81a}, and Item~3.7 below).

   (vii) It is conceivable that many values of $N$ are relevant
for extra-mathematical reasons: 

\par\noindent
   $N \le 12$ pores on pollen-grains (botany); 
\par\noindent
   $N \ge 60$ atoms in various large carbon molecules (chemistry); 
\par\noindent
   $N \sim 20 \, 000$ detectors for a PET tomography of the
brain (medical imaging);
\par\noindent
   $100 \le N \le 10^{20}$ charged particles on a conducting sphere 
(electrostatics).

\bigskip

   In Section 3 below, there are  other examples of cubature formulas and
designs on spheres. 

\bigskip

\noindent {\bf 1.19.\ Several quality criteria for spherical
configurations.} What is the best way to arrange a given number $N$ of
points on a sphere, say here on $\Bbb S^2$ ? 
The answer depends of course of what is meant by \lq\lq the best\rq \rq .
Besides maximizing the strength of the configuration viewed as a spherical
design, we mention here two other criteria.

   The {\bf Tammes' problem} asks what are the configurations
$\{x_1,\hdots,x_n\} \subset \Bbb S^2$ which maxi\-mize
$\min_{1 \le i,j \le n, i\ne j}( \operatorname{distance}(x_i,x_j) )$,
or equivalently minimize
$\max_{1 \le i,j \le n, i\ne j}(\langle x_i \mid x_j \rangle)$.
Configurations of $N$ points on $\Bbb S^2$
are discussed in \cite{Tamme--28} for $N \le 12$.
Pieter L.M.\ Tammes (1903-1980) is a Dutch botanist 
who was interested in the distribution
of places of exit on pollen grains 
(the most frequent case seems to be $N = 3$); 
he should not be confused with 
his aunt Tine Tammes (1871--1947), 
who made important contributions to early genetics \cite{Stamh--95}.

   Here are some of the configurations which are the best from the point
of view of Tammes' problem: 
those of Example~1.18 for $N = 4, 6, 12$
(regular polytopes with triangular faces),
but other configurations for $N = 8$ 
(square antiprisms),
$N = 20$ 
(unknown configurations which are not regular dodecahedras
\cite{vdWae--52}), 
and $N = 24$ (snub cubes, see (\cite{SchWa--51} and \cite{Robin--61}).
More on Tammes' problem in mathematics in \cite{Coxet--62}, 
in Section 35 of \cite{Fejes--64}, 
in Section~2.3 of Chapter 1 of \cite{ConSl--99},
and in Chapter 3 of \cite{EriZi--01}.

   Let $X$ be a finite subset of $\Bbb S^{n-1}$.
When the emphasis is on properties like those of Proposition~1.13,
the set $X$ is called a {\it spherical design}.
However, when the emphasis is on the distance set $D_X$ defined in~1.6, 
the standard name for $X$ is that of {\bf spherical code}.
For many constructions of spherical codes, see \cite{EriZi--01}.

\medskip

   Problem 7 of Smale's {\it Mathematical problems for the next century}
\cite{Smale--98} asks what are the configurations which maximize
$\prod_{1 \le i < j \le N}\operatorname{distance}(x_i,x_j)$.
The problem is motivated by complexity theory, and the search of
algorithms related to the fundamental theorem of algebra.

For a discussion of related criteria, see \cite{SaaKu--97}.

\bigskip

\noindent {\bf 1.20.\ Cubature formulas on other spaces.}
There are documented cubature formulas of strength $3, 5, \hdots$
(with respect to the space of polynomials) on the hypercube
$[-1,1]^n$ of $\Bbb R^n$. 
For compact subsets of the plane with positive area,
there are cubature formulas of strength $k$ 
and of size $\frac{1}{2}(k+1)(k+2)$.
See \cite{DavRa--84}, in particular Section~5.7.

Let $\Omega = \Bbb T^2$ be a $2$-torus of revolution 
embedded in the Euclidean space $\Bbb R^3$, 
together with its standard probability measure $\sigma$
(up to a normalization factor, 
$\sigma$ is the area given by the first fundamental form 
of the surface $\Bbb T^2$ embedded in $3$-space).
Let $\Cal F^{(k)}(\Bbb T^2)$ be the space of functions on the torus
which are restrictions of polynomial functions of total degree 
at most $k$ on $\Bbb R^3$. 
Kuijlaars \cite{Kuijl--95} has shown that there exist constants
$C_1,C_2 > 0$ with the following property: for any $k \ge 0$, there exists
a geometrical design on $\Cal F^{(k)}(\Bbb T^2)$ with a number $N$ of
points satisfying  $C_1 k^2 \le N \le C_2k^2$.

\bigskip
\head
{\bf 
2.\ Cubature formulas and reproducing kernels
}
\endhead
\medskip

   This section begins with a reminder of standard material which goes back
at least to \cite{Arons--50} and \cite{Krein--63}; 
see also \cite{Stewa--76} and \cite{BeHaG--03}.

   Let $\Cal H$ be a real Hilbert space of functions on a set $\Omega$.
We denote by $\langle \varphi_1 \mid \varphi_2 \rangle$ the scalar product
of
$\varphi_1$ and $\varphi_2$ in $\Cal H$. We assume
\footnote{
In all examples appearing below, $\Cal H$ is finite dimensional and this
condition is therefore automatically fulfilled. Hence, the reader can
assume to start with that $\dim_{\Bbb R}(\Cal H) < \infty$.
}
that, for each $\omega \in \Omega$, the evaluation 
$\Cal H \longrightarrow \Bbb R$, $\varphi \longmapsto \varphi(\omega)$
is continuous; this implies that there exists a function 
$\varphi_{\omega} \in \Cal H$ such that 
$\varphi(\omega) = \langle \varphi \mid \varphi_{\omega} \rangle$
for all $\varphi \in \Cal H$.

\bigskip

\proclaim{\bf 2.1.\ Definition} With the notation above, 
the {\bf reproducing kernel} of
$\Cal H$  is the function
$\Phi : \Omega \times \Omega \longrightarrow \Bbb R$ defined by
$\Phi(\omega',\omega) = \langle \varphi_{\omega} \mid
\varphi_{\omega'} \rangle$  
for all $\omega,\omega' \in \Omega$,
or equivalently by $\Phi(\cdot,\omega) = \varphi_{\omega}(\cdot)$.
\endproclaim

The terminology, \lq\lq reproducing\rq\rq , refers to the equality
$\varphi(\omega) = \langle \varphi(\cdot) \mid \Phi(\cdot,\omega)\rangle$
for all $\varphi \in \Cal H$ and $\omega \in \Omega$.

\proclaim{2.2.\ Proposition} Let $\Omega$, $\Cal H$, 
$\left( \varphi_{\omega} \right)_{\omega \in \Omega}$, and $\Phi$ be as
above.

   (i) The kernel $\Phi$ is of positive type. 
In particular, its diagonal values
$\Phi(\omega,\omega)$ are positive, and
$\vert \Phi(\omega,\omega ') \vert ^2 \le 
\Phi(\omega,\omega)\Phi(\omega ',\omega ')$ 
for all $\omega,\omega ' \in \Omega$.

   (ii) The family $\left( \varphi_{\omega} \right)_{\omega \in \Omega}$ 
generates $\Cal H$.

   (iii) If $\Cal H \ne 0$, then $\Phi(\omega,\omega) \ne 0$ for some
$\omega \in \Omega$.

   (iv) If $\left( e_j \right)_{j \in J}$ is any orthonormal basis of
$\Cal H$, then
$$
   \Phi(\omega , \omega ') \, = \, 
   \sum_{j \in J} e_j(\omega)  e_j(\omega ')
$$
for all $\omega,\omega ' \in \Omega$.

   (v) In case there exists a finite positive measure $\sigma$ on
$\Omega$ such that $\Cal H$ is a closed subspace of $L^2(\Omega,\sigma)$,
with $\int_{\Omega \times \Omega} \abs{\Phi(\omega,\omega ')}^2
d\sigma(\omega) d\sigma(\omega ') < \infty$, then
$$
\dim_{\Bbb R}(\Cal H) \, = \,
\int_{\Omega} \Phi(\omega,\omega) d \sigma (\omega) \, < \,  \infty .
$$

  (vi) Assume moreover that $\Omega$ is a topological space and that 
the dimension of $\Cal H$ is finite. Then the following conditions are
equivalent: \par

\hskip.5cm
$(vi_a)$ all functions in $\Cal H$ are continuous;

\hskip.5cm
$(vi_b)$ the kernel 
   $\Phi : \Omega \times \Omega \longrightarrow \Bbb R$ is continuous;

\hskip.5cm
$(vi_c)$ the mapping $\omega \longmapsto \varphi_{\omega}$ 
from $\Omega$ to $\Cal H$ is continuous.

\endproclaim

\demo{Proof} 
$(i)$ By definition, the kernel $\Phi$ is of positive type 
if it is symmetric and if
$$
\sum_{j,k=1}^n \lambda_j \lambda_k \Phi(\omega_j,\omega_k)
\, \ge \, 0
\tag$*$
$$
for all integers $n$, real numbers $\lambda_1,\hdots,\lambda_n$ and 
points $\omega_1,\hdots,\omega_n$ in $\Omega$. 
This is clear here, since the left-hand term 
$
\sum_{j,k=1}^n \lambda_j \lambda_k
\left\langle \varphi_{\omega_k} \mid \varphi_{\omega_j} \right\rangle
$
of $(*)$ is equal to the square of the Hilbert-space norm
of the sum
$
\sum_{k=1}^n  \lambda_k \varphi_{\omega_k}
$.
The last claim of $(i)$ follows by the Cauchy-Schwarz inequality.

   $(ii)$ Observe that, for $\varphi \in \Cal H$, the condition
$\langle \varphi \mid \varphi_\omega \rangle = 0$ 
for all $\omega \in \Omega$ implies that $\varphi(\omega)=0$
for all $\omega \in \Omega$.

   $(iii)$ Assume that $\Phi(\omega,\omega)=0$ for all $\omega \in \Omega$.
Then, by $(i)$, $\Phi(\omega,\omega ')=0$ for all 
$\omega,\omega'\in \Omega$,
and  it follows  from $(ii)$ that $\Cal H = 0$.

   $(iv)$ Evaluate the Fourier expansion
$\varphi_{\omega'} = \sum_{j \in J} 
\langle \varphi_{\omega'} \mid e_j \rangle e_j$
at the point $\omega$, and obtain
$\Phi(\omega,\omega') = 
\sum_{j \in J} \langle \varphi_{\omega'} \mid e_j \rangle
\langle e_j \mid \varphi_{\omega} \rangle
= \sum_{j \in J} e_j(\omega)  e_j(\omega')$.

   $(v)$ Denote by 
$K_{\Phi} : L^2(\Omega,\sigma) \longrightarrow L^2(\Omega,\sigma)$ the
linear operator defined by the kernel $\Phi$, namely by
$(K_{\Phi}\varphi)(\omega) = \int_{\Omega} \Phi(\omega,\omega ')
\varphi(\omega') d\sigma(\omega ')$ for all $\varphi \in \Cal H$ and
$\omega \in \Omega$.  
On the one hand, $K_{\Phi}$ is  a Hilbert-Schmidt operator, 
because of the $L^2$-condition on $\Phi$;
on the other hand, $K_{\Phi}$ is the identity,
since the kernel $\Phi$ is reproducing. 
It follows that the dimension of $\Cal H$ is finite. 
Moreover, we have 
$$
\int_{\Omega} \Phi(\omega,\omega) d\sigma(\omega) \, = \, 
\sum_{j \in J} \int_{\Omega} \abs{e_j(\omega)}^2 d\sigma(\omega) \, = \, 
\sum_{j \in J} \norm{e_j}^2 \, = \, 
\dim_{\Bbb R}(\Cal H)
$$
by $(iv)$.

   $(vi)$ Since $\operatorname{dim}_{\Bbb R}(\Cal H) < \infty$,
the mapping of $(vi_c)$ is continuous if and only if 
the real-valued functions 
$\omega \longmapsto \langle \varphi \mid \varphi_{\omega} \rangle =
\varphi (\omega)$
are continuous for all $\varphi \in \Cal H$,
so that $(vi_a)$ and $(vi_c)$ are equivalent.

  If $(vi_c)$ holds, then $\Phi$ is continuous since it is the composition
of the continuous mapping 
$(\omega ' , \omega) \longmapsto (\varphi_{\omega} , \varphi_{\omega '})$
from $\Omega \times \Omega$ to $\Cal H \times \Cal H$
with the scalar product.
If $(vi_b)$ holds, then $\varphi_{\omega} = \Phi(\cdot,\omega)$
is continuous for all $\omega \in \Omega$,
and $(vi_a)$ follows by $(ii)$.
$\square$ 
\enddemo

\bigskip

\noindent {\bf 2.3.\ Example: atomic measure.} 
Let $X$ be a set and 
$W : X \longrightarrow \Bbb R^*_+$ a positive-valued function. The
Hilbert space $\ell^2 (X,W)$, with scalar product given by
$$
\langle \psi_1 \mid \psi_2 \rangle \, = \, 
\sum_{x \in X} W(x) \psi_1(x)  \psi_2(x)
\quad \text{for all} \quad \psi_1,\psi_2 \in \ell^2(X,W) ,
$$
gives rise to the functions
$$
\psi_{x} = \frac{1}{W(x)}\delta_{x} \, : \, 
y \ \longmapsto \
\left\{
\aligned
W(x)^{-1} \quad &\text{if} \quad y = x \\
0         \qquad &\text{otherwise}
\endaligned
\right.
$$
and to the reproducing kernel $\Psi$ with values
$$
\Psi(x,y) \, = \,
\langle \psi_y \mid \psi_x \rangle \, = \,
\left\{
\aligned
W(x)^{-1} \quad &\text{if} \quad y = x \\
0         \qquad &\text{otherwise.}
\endaligned
\right.
$$
In particular, 
if $X$ is finite and $W(x) = \abs{X}^{-1}$ for all $x \in X$, 
then $\abs{X}^{-1}\Psi$ is the charac\-teristic function of the
diagonal in $X \times X$. 

\bigskip

\noindent {\bf 2.4.\ Example: standard reproducing kernels on spheres.}
Let the notation be as in Example~1.11 and let $k \ge 0$. 
Each of the spaces
$\Cal F^{(k)}(\Bbb S^{n-1})$, $\Cal P^{(k)}(\Bbb S^{n-1})$,
$\Cal H^{(k)}(\Bbb S^{n-1})$
is a subspace of the Hilbert space $L^2(\Bbb S^{n-1},\sigma)$.

   There exists a unique polynomial $Q^{(k)}(T) \in \Bbb R[T]$ 
of degree $k$ with the following pro\-perties:
for any $\omega  \in \Bbb S^{n-1}$, 
the polynomial function defined on $\Bbb R^n$ 
which is homogeneous of degree $k$ 
and of which the restriction to $\Bbb S^{n-1}$ is given by
\footnote{
Where
$\langle \omega \mid \omega ' \rangle = \sum_{i=1}^n
\omega_i \omega'_i$. 
}
$$
\varphi_{\omega } \, : \, 
\omega ' \, \longmapsto \, Q^{(k)}(\langle \omega \mid \omega ' \rangle)
$$
is harmonic, 
and $Q^{(k)}(1) = \dim_{\Bbb R}(\Cal H^{(k)}(\Bbb S^{n-1}))$;
the polynomial $Q^{(k)}$ is a form of a Gegenbauer polynomial
(see, e.g., Theorems~IV.2.12 and~IV.2.14 in \cite{SteWe--71},
or Section~IX.3 in \cite{Vilen--68}). 
It is routine to check that
$Q^{(0)}(X) = 1$, $Q^{(1)}(X) = nX$, and
$$
\frac{k+1}{n+2k}Q^{(k+1)}(X) \, = \,
XQ^{(k)}(X) - \frac{n+k-3}{n+2k-4}Q^{(k-1)}(X)
$$
for $k \ge 2$.
Observe that $Q^{(k)}(T)$ is of the form 
$$
Q^{(k)}(T) \, = \, \sum_{j=0}^{[k/2]} (-1)^j c_j^{(k)} T^{k-2j}
$$
with $c_0^{(k)}, c_1^{(k)}, \hdots > 0$.
The reproducing kernel of the space $\Cal H^{(k)}(\Bbb S^{n-1})$ is given
by
$$
\Phi^{(k)}(\omega,\omega ') \, = \, 
Q^{(k)}(\langle \omega \mid \omega ' \rangle)
$$
for all $\omega,\omega ' \in \Bbb S^{n-1}$.
It is a restriction of the kernel
$$
\Bbb R^n \times \Bbb R^n  \ni  (\omega,\omega ') 
\, \longmapsto \, 
\sum_{j=0}^{[k/2]} (-1)^j c_j^{(k)} 
\langle \omega   \mid \omega ' \rangle ^{k-2j}
\langle \omega   \mid \omega   \rangle ^{j}
\langle \omega ' \mid \omega ' \rangle ^{j}
\, \in \, \Bbb R
$$
which is homogeneous of degree $k$ in each variable $\omega,\omega '$
separately.

   The reproducing kernels of the spaces $\Cal P^{(k)}(\Bbb S^{n-1})$ and 
$\Cal F^{(k)}(\Bbb S^{n-1})$ are given similarly in terms of the
polynomials
$$
C^{(k)}(T) \, = \, \sum_{j=0}^{[k/2]} Q^{(k-2j)}(T) 
\qquad \text{and} \qquad
R^{(k)}(T) \, = \, \sum_{j=0}^{k} Q^{(j)}(T) .
$$
Observe that
$$
\aligned
\dim_{\Bbb R}(\Cal H^{(k)}(\Bbb S^{n-1})) \, &= \, Q^{(k)}(1) \, = \, 
\binom{n+k-1}{n-1} - \binom{n+k-3}{n-1} \\
\dim_{\Bbb R}(\Cal P^{(k)}(\Bbb S^{n-1})) \, &= \, C^{(k)}(1) \, = \, 
\binom{n+k-1}{n-1}  \\
\dim_{\Bbb R}(\Cal F^{(k)}(\Bbb S^{n-1})) \, &= \, R^{(k)}(1) \, = \, 
\binom{n+k-1}{n-1} + \binom{n+k-2}{n-1} 
\endaligned
$$
by Proposition~2.2.v.

   This material is standard: see, e.g., \cite{SteWe--71} and
\cite{DeGoS--77}.

\bigskip

   Proposition~2.6 is the first general existence result 
for cubature formulas of this exposition. 
(It is a minor strengthening of 
Theorems~2.8 and~3.17 in \cite{LyuVa--93}.)

\bigskip

\proclaim{2.5.\ Lemma} 
Let $V$ be a finite dimensional real vector space,
let $\mu$ be a probability measure on $V$,
and let $b = \int_V vd\mu(v)$ denote its barycentre.
If $C$ is a convex subset of $V$ such that $\mu(C) = 1$,
then $b \in C$.
\endproclaim

\demo{Remarks} (i) The lemma is probably well-known, but we haven't been
able to trace it in print. The proof below was shown to us by Yves
Benoist.

(ii) The point of the lemma is that $b$ is not only in the closure of $C$,
but in $C$ itself.
\enddemo

\demo{Proof} Let $U$ be the minimal affine subspace of $V$ such that
$\mu(C \cap U) = 1$; upon restricting $C$, 
we can assume that $C \subset U$. 
Upon translating $\mu$ and $C$,
we can assume that $b = 0$.

   We claim that $b = 0$ is in the interior of $C$ inside $U$.
If this was not true, there would exist by the Hahn-Banach theorem
a non-zero linear form $\xi$ on $U$ 
such that $\xi (c) \ge 0$ for all $c \in C$.
The equality
$$
0 \, = \, \xi (b) \, = \, 
\int_{U} \xi(w) d\mu(w) \, = \, \int_{C} \xi(c) d\mu(c) 
$$
would imply $\xi(c) = 0$ for almost every $c \in C$;
thus, we would have 
$\mu(C \cap \operatorname{Ker}(\xi)) = 1$, in contradiction with the
definition of $U$.

   Hence $b$ is in the interior of $C$ inside $U$, 
and in particular $b \in C$.
$\square$   
\enddemo

\bigskip

\proclaim{2.6.\ Proposition} 
Let $(\Omega,\sigma)$ be a finite measure space. 
Let $\Cal H$ be a finite dimensional Hilbert space of functions on
$\Omega$ which is a subspace of $L^2(\Omega,\sigma)$
and which contains the constant functions.

   Then there exists a cubature formula $(X,W)$ for $\Cal H$ such that
$$
\abs{X} \, \le \, \dim_{\Bbb R}(\Cal H) .
$$
\endproclaim

\demo{Proof}
Assume first that $\sigma$ is a probability measure.
Let $\Cal H^0$ be the orthogonal supplement of the constants in
$\Cal H$ and let $\Phi^0$ denote its reproducing kernel.
Recall that, for $\omega \in \Omega$, the function
$\varphi_{\omega}^0 : \omega ' \longrightarrow 
\Phi^0(\omega ',\omega)$ is in $\Cal H^0$.
The set
$$
\tilde{\Omega} \, = \,
\left\{ \varphi \in \Cal H^0 \mid \varphi = \varphi^0_{\omega}
\quad \text{for some} \quad \omega \in \Omega \right\}
$$
linearly generates $\Cal H^0$ by Proposition~2.2.ii.
Observe that
$$
\int_{\Omega} \varphi_{\omega}^0(\omega') d\sigma(\omega)
\, = \, 
\int_{\Omega}  \varphi_{\omega'}^0(\omega) d\sigma(\omega)
\, = \, 
\langle 1 \mid \varphi_{\omega'}^0 \rangle_{\Cal H^0}
\, = \, 
0
$$
for all $\omega ' \in \Omega$,
by definition of $\Cal H^0$. 
It follows from the previous lemma,
applied to the image of $\sigma$ on $\tilde \Omega$,
that $0$ is in the convex hull of $\tilde{\Omega}$.

By Carath\'eodory's theorem 
(see~e.g.~Theorem~11.1.8.6 in~\cite{Berge--78}), 
there exist a finite subset $X$ of $\Omega$ of cardinality at most
$\dim_{\Bbb R}(\Cal H ^0) +1 = \dim_{\Bbb R}(\Cal H)$ 
and a weight function $W : X \longrightarrow \Bbb R_+^*$ such that
$$
\sum_{x \in X} W(x) \, = \, 1 \qquad \text{and} \qquad
\sum_{x \in X} W(x) \varphi_x^0 \, = \, 0 \, \in \, \Cal H^0 .
$$
For all $x \in X$ and $\omega ' \in \Omega$, we have
$\varphi_x^0(\omega ') = 
\Phi^0(\omega ' , x) =
\Phi^0(x, \omega ') = 
\varphi^0_{\omega '}(x)$; hence
$$
\sum_{x \in X} W(x) \varphi_{\omega '}^0(x) 
\, = \, 0 \, = \, 
\int_{\Omega} \varphi_{\omega '}^0 (\omega) d\sigma(\omega) 
$$
for all $\omega ' \in \Omega$.
Since $\left( \varphi_{\omega '}^0 \right)_{\omega ' \in \Omega}$
generates $\Cal H^0$, we have
$$
\sum_{x \in X} W(x) \varphi (x) 
\, = \, 
\int_{\Omega} \varphi (\omega) d\sigma(\omega)
$$
for all $\varphi \in \Cal H^0$. 
Since the same equality holds obviously for constant functions,
it holds also for all $\varphi \in \Cal H$.

In case $\sigma$ is not a probability measure, we can multiply
all values of $W$ by $\sigma(\Omega)$.
Then the previous equality holds again for constant functions on $\Omega$,
and a posteriori for all $\varphi \in \Cal H$.
$\square$
\enddemo

In the situation of Proposition 2.6, 
suppose moreover that we have a cubaure formula $(X',W')$ for $\Cal H$.
Then there  exists a cubature formula $(X,W)$ for $\Cal H$ such that
$$
X \subset X'
\qquad \text{and}Ê\qquad
\abs{X} \, \le \, \dim_{\Bbb R}(\Cal H) .
$$
(This follows from the proof in \cite{Berge--78} of Carath\'eodory's theorem.
Alternatively, we can apply Proposition 2.6 to the subspace of $\ell^2(X',W')$
of restrictions to $X'$ of functions in $\Cal H$.)

\bigskip

\proclaim{2.7.\ Proposition} In the situation of the previous proposition,
assume moreover that $\Omega$ is a connected topological space
and  that $\Cal H$ is a space of continuous functions.

Then the bound on the size of $X$ can be improved to
$$
\abs{X} \, \le \, \dim_{\Bbb R}(\Cal H) - 1 .
$$
\endproclaim

\demo{Proof} 
With the notation of the proof of Proposition~2.6,
the subspace $\tilde \Omega$ of $\Cal H ^0$ is a continuous image
of $\Omega$ by Proposition~2.2.vi, 
so that $\tilde \Omega$ is connected.
In this case, the Carath\'eodory bound on $X$ can be lowered by $1$, 
by a classical theorem of Fenchel. 
See, e.g., Theorem 18 in \cite{Eggle--58}.
$\square$
\enddemo

\bigskip

The following result is an immediate consequence of Propositions 1.7 and 2.7.

\bigskip

\proclaim{2.8.\ 
Theorem (existence of cubature formulas, with bounds on sizes)} 
Let $\Omega$ be a connected topological space 
with a finite  measure $\sigma$ and
let $\left(\Cal F^{(k)} \right)_{k \ge 0}$ be 
a sequence of polynomial spaces on $(\Omega,\sigma)$.
Assume that functions in $\bigcup_{k \ge 0} \Cal F^{(k)}$
are continuous on $\Omega$. Choose $l \ge 0$.

Then there exists a finite subset $X$ of $\Omega$ 
such that
$$
\dim_{\Bbb R}(\Cal F^{(l)}) \, \le \, 
\abs{X} \, \le \, 
\dim_{\Bbb R}(\Cal F^{(2l)}) - 1.
$$
and a weight $W : X \longrightarrow \Bbb R_+^*$
such that $(X,W)$ is a cubature formula of strength $2l$
on $\Omega$.
\endproclaim

\bigskip

   Our next proposition (2.10) shows properties 
of cubature formulas $(X,W)$
for which the previous lower bound is an equality, 
namely which are tight (Definition~1.8).

\bigskip

\proclaim{2.9.\ Lemma} Consider  data consisting of
\roster
\item"(a)" 
a finite measure space $(\Omega,\sigma)$, 
a Hilbert space of functions
$\Cal H \subset L^2(\Omega,\sigma)$, 
and the corresponding reproducing kernel $\Phi$;
\item"(b)" 
a finite subset $X$ of $\Omega$, 
a weight $W : X \longrightarrow \Bbb R_+^*$,
the Hilbert space   $\ell^2(X,W)$, 
and the kernel $\Psi$ as in Example~2.3;
\item"(c)" 
the restriction mapping 
$\rho : \Cal H \longrightarrow \ell^2(X,W)$,
$\varphi \longmapsto \varphi \mid_X$,
and the adjoint mapping 
$\rho^* : \ell^2(X,W) \longrightarrow \Cal H$.
\endroster
With the previous notation, namely with $\varphi_{\omega}$
as in 2.1 and $\psi_x$ as in 2.3, we have
\roster
\item"(i)" 
   $\rho^*(\psi_x) = \varphi_x$ for all $x \in X$;
\item"(ii)" 
   $\rho^*$ is an isometry if and only if
   $\Psi$ is the restriction of $\Phi$ to $X \times X$;
\item"(iii)" 
   in case $\dim_{\Bbb R}(\Cal H) = \abs{X}$, the mapping $\rho$
   is an isometry if and only if
   $\Psi$ is the restriction of $\Phi$ to $X \times X$.
\endroster
\endproclaim

\demo{Proof} (i) By definition of $\rho^*$, we have
for all $\varphi \in \Cal H$ and for all $x \in X$
$$
\langle \varphi \mid \rho^*(\psi_x) \rangle \, = \, 
\langle \rho(\varphi) \mid \psi_x \rangle \, = \,
(\rho(\varphi))(x) \, = \, \varphi(x)
$$
and therefore $\rho^*(\psi_x) = \varphi_x$.

   (ii) Let $x,y \in X$. On the one hand,
$\langle \rho^*(\psi_x) \mid \rho^*(\psi_y) \rangle = \Phi(x,y)$
by (i); on the other hand,
$\langle \psi_x \mid \psi_y \rangle = \Psi(x,y)$.
Claim (ii) follows.

   (iii) A linear mapping between two finite-dimensional Hilbert spaces
of the same dimension is an isometry if and only if its adjoint is an
isometry~; thus (iii) follows from (ii).
$\square$

\enddemo

\bigskip

\proclaim{2.10.\ Proposition}  
Let $(\Omega,\sigma)$ be a finite measure space and let
$( \Cal F^{(k)} )_{k \ge 0}$ be a sequence of polynomial spaces on
$(\Omega,\sigma)$;
assume that the corresponding reproducing kernels $\Phi^{(k)}$
are constant on the diagonal of $\Omega$.

Choose $l\ge 0$ and let $(X,W)$ be a cubature formula of strength $2l$
on $\Omega$ which is tight. 

Then $W$ is uniform, namely $X$ is a tight
geometrical $(2l)$-design. 
\endproclaim

\demo{Proof} The restriction mapping
$$
\rho : \Cal F^{(l)} \longrightarrow \ell^2(X,W)
$$
is an isometry by Proposition~1.7. 
If $(X,W)$ is tight, 
the domain and the range of $\rho$ have the same dimension, 
so that $\rho$ is onto. 
Thus, if $\Psi : X \times X \longrightarrow \Bbb R$
denotes as in Example~2.3 the reproducing kernel of $\ell^2(X,W)$, 
then
$$
\Psi(x,x') \, = \, \Phi^{(l)}(x,x') 
\qquad \text{for all} \quad x,x' \in X
$$
by Lemma~2.9.
It follows that $\Psi$ is constant on the diagonal of 
$X \times X$, so that the weight $W$ is constant by Example~2.3; 
otherwise said, $X$ is a geometrical $(2l)$-design.
$\square$
\enddemo

\bigskip

\proclaim{2.11.\ Definition} 
Let $\Omega$ be a metric space, 
let $\sigma$ be a finite measure on $\Omega$, 
let $\Cal H$ be a finite dimensional Hilbert space of functions on
$\Omega$ which is a subspace of $L^2(\Omega,\sigma)$, 
and let $\Phi$ denote the corresponding reproducing kernel.
We say that {\bf Condition (M) holds for $\Cal H$ 
and a function $S$}
if the values  of the  kernel 
depend only on the distances:
$$
\Phi(\omega,\omega ') \, = \, S(d(\omega,\omega '))
$$
for all $\omega,\omega ' \in \Omega$, and for some function
$S : \Bbb R_+ \longrightarrow \Bbb R$. 
\endproclaim

   Observe that, by Proposition~2.2.v, this implies that
$$
\Phi(\omega,\omega) \, = \, S(0) \, = \, 
\sigma(\Omega)^{-1} \dim_{\Bbb R}(\Cal H)
$$ 
for all $\omega \in \Omega$.

   On spheres, Condition (M) holds
\footnote{
Or rather the analogous condition for scalar products
$\langle \omega \mid \omega ' \rangle$
rather than distances $d(\omega,\omega ')$.
} 
for the spaces 
$\Cal H^{(l)}(\Bbb S^{n-1})$ [respectively
$\Cal P^{(l)}(\Bbb S^{n-1})$, $\Cal F^{(l)}(\Bbb S^{n-1})$]
with  $S = Q^{(l)}$ [respectively $S = C^{(l)}$, $S = R^{(l)}$];
the notation is that of Example~2.4.

\bigskip

\proclaim{2.12.\ Proposition} Let $(\Omega,\sigma)$, 
$( \Cal F^{(k)} )_{k \ge 0}$, and $(\Phi^{(k)})_{k \ge 0}$
be as in Proposition~2.10.
Choo\-se $l \ge 0$ and assume that  
Condition (M) holds for $\Cal F^{(l)}$
and a function $S^{(l)}$.

   If $X$ is a tight geometrical $2l$-design on $\Omega$, 
the set of distances
$$
D_X \, = \, \left\{ c \in \mathopen]0,\infty\mathclose[ \ \big\vert  \
c = d(x,y) \quad \text{for some} \quad x,y \in X, x \ne y \right\}
$$
is contained in the set of zeros of the function $S^{(l)}$.
\endproclaim

\demo{Proof}
Notation being as in the proof of Proposition~2.10,
we have 
$$
\Psi(x,x') \, = \, \Phi^{(l)}(x,x') \, = \,
S^{(l)}(d(x,x'))
$$
for all $x,x' \in X$.
Since $\Psi(x,x') = 0$ when $x \ne x'$ (see Example~2.3),
the proposition follows. 
$\square$
\enddemo

\bigskip

   The last proposition of this section shows how the setting of
reproducing kernels makes it straightforward to generalize a
characterization which is well-known for spherical designs
(Theorem~5.5 in \cite{DeGoS--77}, 
or Theorems~3.2 and~4.3 in \cite{GoeSe--79}).

\proclaim{2.13.\ Proposition} 
Let $(\Omega,\sigma)$ be a finite measure space, 
let $\Cal H$ be a finite dimensional Hilbert space of functions on
$\Omega$ which is a subspace of $L^2(\Omega,\sigma)$, 
and let $\Phi$ denote the corresponding reproducing kernel.
Let $X$ be a non-empty finite subset of $\Omega$ and let $W$
be a weight function on $X$.

   (i) We have $\sum_{x,y \in X} W(x)W(y)\Phi(x,y) \, \ge \, 0$.
\smallskip

   (ii)  Equality holds in (i) if and only if
$\sum_{x \in X} W(x) \varphi(x) = 0$ for all $\varphi \in \Cal H$.
\endproclaim

\demo{Proof}
The inequality holds in (i) since $\Phi$ is a kernel of positive type.
 
Assume now equality. For any orthonormal basis $(e_j)_{j \in J}$ of   
$\Cal H$, Proposition~2.2.iv implies that
$$
\aligned
\sum_{x,y \in X} W(x)W(y)\Phi(x,y) \, &= \,
\sum_{x,y \in X} W(x)W(y) \sum_{j \in J} e_j(x)e_j(y) \, 
\\
&= \, \sum_{j \in J} \Big( \sum_{x \in X} W(x)e_j(x) \Big)^2 \, = \, 0 .
\endaligned
$$
Hence $\sum_{x \in X} W(x) e_j(x) = 0$ for all $j \in J$, 
and therefore
$\sum_{x \in X} W(x) \varphi(x)  =  0$
for all $\varphi \in \Cal H$.

The converse implication in (ii) is straightforward.
$\square$
\enddemo

\proclaim{2.14.\ Particular case}
Consider two integers $n \ge 2$ and $k \ge 1$, 
a non-empty finite subset $X$ of $\Bbb S^{n-1}$,
a weight $W : X \longrightarrow \Bbb R_+^*$ such that 
$\sum_{x \in X}W(x) = 1$,
and the polynomials $Q^{(j)}$, $R^{(j)}$ defined in Example 2.4.
Then
$$
\sum_{x,y \in X} W(x)W(y) Q^{(j)}(\langle x \mid y \rangle) 
\, \overset{(*)}\to{\ge} \, 0
\hskip.5cm \text{and} \hskip.5cm
\sum_{x,y \in X} W(x)W(y) R^{(j)}(\langle x \mid y \rangle)
\, \overset{(**)}\to{\ge} \, 1
$$
for all $j \ge 0$. Moreover, the following properties are equivalent: \par
(i) $(X,W)$ is a cubature formula of strength $k$ on $\Bbb S^{n-1}$; \par
(ii) equality holds in $(*)$ for all $j \in \{1,\hdots,k\}$; \par
(iii) equality holds in $(**)$ for $j = k$.
\endproclaim

\bigskip
\head
{\bf 
3.\ The case of spheres \\ Antipodal cubature formulas and spherical
designs
}
\endhead
\medskip

   In this section, we assume that $\Omega = \Bbb S^{n-1}$ is the unit
sphere in $\Bbb R^n$, $n \ge 2$, that $\sigma$ is the rotation-invariant
probability measure on $\Bbb S^{n-1}$, and that the notation is as in
Examples~1.11 and~2.4.

   For tight $(2l)$-spherical designs, 
the set $D_X$ of Proposition~2.12 is precisely known, 
and $\abs{D_X} = l$. 
Instead of distances and $D_X$,
it is more convenient to use scalar products and the set
$$
A_X \, = \, 
\left\{ c \in [-1,1[ \ \big\vert  \
c = \langle x \mid y \rangle 
\quad \text{for some} \quad 
x,y \in X, \ x \ne y \right\} .
$$
The following proposition is Theorem~5.1 of \cite{GoeSe--79}.

\bigskip

\proclaim{3.1.\ Proposition} If $X$ is a tight spherical $(2l)$-design on
$\Bbb S^{n-1}$, 
then $A_X$ coincides with the set of roots of the polynomial
$R^{(l)}$ defined in~2.4.  
\endproclaim

\demo{Proof} Let $Z^{(l)}_R$ denote the set of roots of $R^{(l)}$, 
which is of order $l$. 
Since $A_X \subset Z^{(l)}_R$ by the proof of Proposition~2.12, 
it is enough to show that $a = \abs{A_X}$ is not less than $l$.
More generally, let us show that $a \ge l$ for any spherical
$(2l)$-design.

   Define for each $x \in X$ a polynomial function $\gamma_x$
of degree $a$ by
$$
\gamma_x(\omega) \, = \, \prod_{c \in A_X} 
\frac{\langle x \mid \omega \rangle - c}{1-c} \quad
\text{for all} \quad \omega \in \Bbb R^n .
$$
Then $\gamma_x(x') = \delta_{x,x'}$ for all $x,x' \in X$, so that the 
family  $\left( \gamma_x \right)_{x \in X}$ is linearly independent. Thus
$$
\dim_{\Bbb R}(\Cal F^{(a)}(\Bbb S^{n-1})) \, \ge \, \abs{X} .
$$
Since $\abs{X} \ge \dim_{\Bbb R}(\Cal F^{(l)}(\Bbb S^{n-1}))$ 
by Proposition~1.7,  we have $a  \ge l$. 

In particular, $A_X = Z^{(l)}_R$ for a tight spherical $(2l)$-design.
$\square$
\enddemo

\bigskip

   Propositions~1.7, 2.10, 2.12, 3.1, and Theorem~2.8 apply to 
cubature formulas and spherical designs of {\it even strength}. 
We expose now the analogous facts for {\it odd strength}.

\bigskip

\proclaim{3.2.\ Proposition} {\rm (Compare with~1.7.)} 
Let $n \ge 2$, $l \ge 0$,  
and $\Cal P^{(l)}(\Bbb S^{n-1}) \subset L^2(\Bbb S^{n-1},\sigma)$ be
as in Example~1.11. Let moreover $X = Y \sqcup (-Y)$ be a non-empty
antipodal finite subset of the sphere and  
$W : X \longrightarrow \Bbb R_+^*$ be a symmetric weight. 

Then $(X,W)$ is a cubature formula of strength $2l+1$
if and only if the restriction mapping 
$$
\rho \, : \, \Cal P^{(l)} \ \longrightarrow \ \ell^2(Y,2W)
$$
is an isometry. 
In particular, if $(X,W)$ is a cubature formula of strength $2l+1$, then
$$
\abs{X} \, \ge 2 \, \dim_{\Bbb R}(\Cal P^{(l)}(\Bbb S^{n-1})) \, = \, 
2\binom{n+l-1}{n-1} . \tag{$\sharp$}
$$ 
\endproclaim

\demo{Proof} Observe that 
$\psi(-\omega) = (-1)^j\psi(\omega)$
for all $j \ge 0$, $\psi \in \Cal P^{(j)}(\Bbb S^{n-1})$, 
and $\omega \in \Bbb S^{n-1}$,
so that
$$
\sum_{x \in X} W(x) \varphi_1(x)  \varphi_2(x)
\, = \, 
2\sum_{y \in Y} W(y) \varphi_1(y)  \varphi_2(y)
$$
for $\varphi_1,\varphi_2 \in \Cal P^{(l)}(\Bbb S^{n-1})$.

   The condition for $\rho$ to be an isometry is
$$
2\sum_{y \in Y} W(y) \varphi_1(y)  \varphi_2(y)
\, = \, 
\int_{\Bbb S^{n-1}} \varphi_1(\omega)  \varphi_2(\omega)
d\sigma(\omega)
$$
for all $\varphi_1,\varphi_2 \in \Cal P^{(l)}(\Bbb S^{n-1})$. 
As in the proof of Proposition~1.7, this can be written
$$
2 \sum_{y \in  Y} W(y) \varphi(y)
\, = \, 
\int_{\Bbb S^{n-1}} \varphi(\omega) d\sigma(\omega)
$$
or indeed
$$
\sum_{x \in X} W(x) \varphi(x)
\, = \, 
\int_{\Bbb S^{n-1}} \varphi(\omega) d\sigma(\omega) \tag$*$
$$
for all $\varphi \in \Cal P^{(2l)}(\Bbb S^{n-1})$. 
Now $(*)$ holds for all $\varphi \in \Cal P^{(2l+1)}(\Bbb S^{n-1})$,
since both terms are zero in this case by symmetry reasons.
Hence $\rho$ is an isometry if and only if
$(*)$ holds for all $\varphi \in \Cal F^{(2l+1)}(\Bbb S^{n-1})$.
$\square$
\enddemo

\bigskip
\noindent
{\bf 3.3.\ Remarks and definition.} (i) Let $(X,W)$ be an antipodal
cubature formula  on $\Bbb S^{n-1}$, and write $X~= Y~\sqcup~(-Y)$.  
Observe that, since a function on $X$ can be canonically written 
as the sum of an even function and an odd function,
we have a decomposition of $\ell^2(X,W)$ as the orthogonal direct sum of
two copies of $\ell^2(Y,2W)$. 
It follows from Propositions~1.7 and~3.3 
that the restriction mapping
$$
   \Cal F^{(l)}(\Bbb S^{n-1}) 
   \, = \,  
   \Cal P^{(l)}(\Bbb S^{n-1}) \oplus \Cal P^{(l-1)}(\Bbb S^{n-1}) \
\longrightarrow \ \ell^2(X,W) = \ell^2(Y,2W) \oplus \ell^2(Y,2W)
$$
is an isometry if and only if the restriction mapping
$
\Cal P^{(l)}(\Bbb S^{n-1}) \longrightarrow  \ell^2(Y,2W)
$
is an iso\-metry. 

\medskip

   (ii) For a spherical $(2l+1)$-design $X$ in $\Bbb S^{n-1}$
(not necessarily antipodal), it is known that
$\abs{X} \ge 2 \dim_{\Bbb R}(\Cal P^{(l)}(\Bbb S^{n-1}))$;
moreover, in case equality holds, then $X$ is antipodal.
See Theorem~5.12 in \cite{DeGoS--77};
the proof uses the \lq\lq linear programming method\rq\rq ,
which is more powerful than the \lq\lq Fisher type method\rq\rq
\ used in our proof of Proposition~3.2. 

\medskip 

   (iii)  A cubature formula $(X,W)$ of strength
$2l+1$ on $\Bbb S^{n-1}$ is {\bf tight} (see the definition in~1.11)
if  equality holds in Equation $(\sharp)$ of Proposition~3.2.

\bigskip

The following is Theorem~5.11 of \cite{GoeSe--79}.

\proclaim{3.4.\ Proposition} {\rm (Compare with~2.10, 2.12, and 3.1.)} 
Let $(X,W)$ be a cubature formula of strength $2l+1$ on $\Bbb S^{n-1}$
which is tight.

   Then $W$ is uniform, namely $X$ is a tight spherical $(2l+1)$-design.
Moreover the set
$$
B_X \, = \, 
\left\{ c \in \mathopen]-1,1\mathclose[ \ \big\vert  \
c = \langle x \mid y \rangle 
\quad \text{for some} \quad 
x,y \in X, \ x \ne \pm y \right\} 
$$
coincides with the set of roots of the polynomial $C^{(l)}$ defined in~2.4.
\endproclaim

\demo{Proof} Let $Z^{(l)}_C$ denote the set of roots of $C^{(l)}$,
which is of order $l$.
Since $B_X \subset Z^{(l)}_C$ by the argument used in the proof of
Proposition~2.12, it is enough to show that $b = \abs{B_X}$ 
is not less than $l$.
More generally, let us show that $b \ge l$ for any
antipodal $(2l+1)$-design.

Define for each $x \in X$ a polynomial function 
$\tilde \gamma_x$ of degree $b$ by
$$
\tilde \gamma_x(\omega) \, = \, \prod_{c \in B_X} 
\frac{\langle x \mid \omega \rangle - c}{1-c} \quad
\text{for all} \quad \omega \in \Bbb R^n .
$$
As $-B_X = B_X$, we have
$$
\tilde \gamma_x(\omega) \, = \, 
\frac{ \prod_{c \in B_X} (\langle x \mid \omega \rangle + c) }
{ \prod_{c \in B_X} (1-c) } \, = \, 
(-1)^b
\frac{ \prod_{c \in B_X} (\langle x \mid -\omega \rangle - c) }
{ \prod_{c \in B_X} (1-c) } \, = \, 
(-1)^b \tilde \gamma_x (-\omega)
$$
for all $\omega \in \Bbb R^n$. 
Thus the restriction of $\tilde \gamma_x$ to $\Bbb S^{n-1}$ 
is in $\Cal P^{(b)}(\Bbb S^{n-1})$.
Decompose $X$ as $Y \sqcup (-Y)$ as in Item~1.12;
since $\tilde \gamma_y(y') =  \delta_{y,y'}$
for all $y,y' \in Y$, we have
$$
\operatorname{dim}_{\Bbb R} \big( \Cal P^{(b)}(\Bbb S^{n-1}) \big) 
\, \ge \,  \abs{Y} \, \ge \,
\operatorname{dim}_{\Bbb R} \big( \Cal P^{(l)}(\Bbb S^{n-1}) \big)
$$
and therefore $b \ge l$.

   In particular $B_X = Z^{(l)}_C$ for a tight antipodal spherical
$(2l+1)$-design.
$\square$
\enddemo

\bigskip
\noindent
{\bf 3.5.\ Remark.} Given $n$ and $l$,
let $(X,W)$ be a cubature formula of strength $2l+1$ on $\Bbb S^{n-1}$
such that $\abs{X} \le \abs{X'}$
for any cubature formula of strength $2l+1$ on $\Bbb S^{n-1}$.
The weight $W$ need not be uniform;
see Example~3.7 for such a cubature formula
of strength $7$ on $\Bbb S^2$.

\bigskip

\proclaim{\bf 3.6.\ Proposition} {\rm (Compare with~2.8.)}
For each $l \ge 0$, there exists a finite antipodal subset
$X = Y \sqcup (-Y)$ of $\Bbb S^{n-1}$ such that
$$
\operatorname{dim}_{\Bbb R} \big( P^{(l)}(\Bbb S^{n-1}) \big) 
\, \le \, \abs{Y} \, \le \,
\operatorname{dim}_{\Bbb R} \big( P^{(2l)}(\Bbb S^{n-1}) \big) - 1
$$
and a symmetric weight $W : X \longrightarrow \Bbb R_+^*$
such that $(X,W)$ is an antipodal cubature formula of
strength $2l+1$ on $\Bbb S^{n-1}$.
\endproclaim

\demo\nofrills{Proof: \usualspace} see Propositions~2.6 and 2.7.
$\square$

Similar estimates are known for other spaces, in particular for hypecubes
in $\Bbb R^n$; see page 366 of \cite{DavRa--84}.
\enddemo

\bigskip

\noindent {\bf 3.7.\ Example: cubature formula of strength $5$, $7$, and $9$
on  $\Bbb S^2$.}  
A number of explicit cubature formulas can be collected from the
litera\-ture, either directly or indirectly. Many of them are reviewed
in \cite{Dicks--19} (see in particular pages 717--724) and
\cite{Rezni--92}. Let us first write down an identity of Lucas (1877):
$$
8(u_1^4 + u_2^4 + u_3^4) + 
\sum_{\epsilon_2,\epsilon_3 \in \{1,-1\}}
   (u_1 + \epsilon_2u_2 + \epsilon_3u_3)^4 \, = \,
12(u_1^2 + u_2^2 + u_3^2)^2 .
$$
Let $Y \subset \Bbb S^2$ be the set of size $7$ 
containing the three vectors $e_1, e_2, e_3$ 
of the canonical orthonormal basis of $\Bbb R^3$
and the four vectors
$3^{-1/2}(e_1 + \epsilon_2e_2 + \epsilon_3e_3)$,
for $\epsilon_2,\epsilon_3 = 1,-1$.
Define a weight $W : Y \longrightarrow \Bbb R_+^*$
by $W(e_i) = \frac{8}{60}$ and
$W(3^{-1/2}(e_1 + \epsilon_2e_2 + \epsilon_3e_3))= \frac{9}{60}$.
Then Lucas' identity can be rewritten as
$$
\sum_{y \in Y} W(y) \langle y \mid u \rangle ^4 
\, = \,
\frac{1}{5} \langle u \mid u \rangle ^2
\quad \text{for all} \quad u \in \Bbb R^3
$$
so that $(Y,W)$ is a cubature formula of size $7$ 
for $\Cal P^{(4)}(\Bbb S^2)$
by the argument of Proposition~1.13.
If $X = Y \sqcup (-Y)$ and if $W$ is extended by symmetry, 
$(X,\frac{1}{2}W)$ is
a cubature formula of size $14$ and strength $5$.

\smallskip

   Size $\abs{Y} = 7$ [respectively $\abs{X} = 14$] is not optimal for 
cubature formulas for $\Cal P^{(4)}(\Bbb S^2)$ 
[respectively $\Cal F^{(5)}(\Bbb S^2)$].
Indeed, as already reported in Items~1.17 and 1.18,  
the $12$ vertices of a regular icosahedron provide a tight spherical
$5$-design, 
and such a design is unique up to isometry 
\cite{DeGoS--77, page 375}.
The $6$ vertices of a corresponding $Y$ provide a design
for the space $\Cal P^{(4)}(\Bbb S^2)$. 
(Note that this is not a spherical $4$-design, 
because $\frac{1}{\abs{Y}} \sum_{y \in Y}\varphi(y)$
is not always equal to $\int_{\Bbb S^2} \varphi(\omega) d\sigma(\omega)$
for functions $\varphi$ in the second summand of the decomposition 
$\Cal F^{(4)}(\Bbb S^2) =
\Cal P^{(4)}(\Bbb S^2) \bigoplus \Cal P^{(3)}(\Bbb S^2)$.)

\smallskip

From \cite{Rezni--92} and \cite{HarSl--96} (see also Item~1.18), 
we collect the following facts.

   $\circ$ There exists a pair $(Y,W)$ with $\abs{Y} = 11$ 
which is a cubature formula for $\Cal P^{(6)}(\Bbb S^2)$,
and $11$ is the minimal possible size 
(pages~133--135 of \cite{Rezni--92}).
Hence there exists a cubature formula $(X,\frac{1}{2}W)$ for 
$\Cal F^{(7)}(\Bbb S^2)$ of size $\abs{X} = 22$.
On the other hand, a spherical design $X \subset \Bbb S^2$
of size $\abs{X} \le 22$ is of strength at most $5$,
and a spherical $7$-design in $\Bbb S^2$ has size at least $24$
(compare with the lower bound of Proposition~3.3 
for $7$-designs in $\Bbb S^2$,  which is~$20$).

   $\circ$ There exists a pair $(Y,W)$ with $\abs{Y} = 16$ 
which is a cubature formula for $\Cal P^{(8)}(\Bbb S^2)$,
and $16$ is the minimal possible size 
(pages~111 and~136 of \cite{Rezni--92}, 
referring to Finden, Sobolev, and McLaren).
Hence there exists a cubature formula $(X,\frac{1}{2}W)$ for 
$\Cal F^{(9)}(\Bbb S^2)$ of size $\abs{X} = 32$. 
On the other hand, a spherical design $X \subset \Bbb S^2$ 
of size $\abs{X} \le 32$ is of strength at most $7$,
and a spherical $9$-design in $\Bbb S^2$ has size at least $48$
(compare with the lower bound of Proposition~3.2 
for $9$-designs in $\Bbb S^2$,   which is $30$).

\bigskip

\noindent {\bf 3.8.\ Example: cubature formulas of strengths $5$, $7$,
$9$ and $11$ on  $\Bbb S^3$.}    
The following identities go back respectively to
Liouville (1859), Kempner (1912), Hurwitz (1908), and J. Schur (1909):
$$
\sum_{i=1}^4 (2u_i)^4 \, + \,
\sum_{\epsilon_i \in \{1,-1\}}
   (u_1 + \epsilon_2 u_2 + \epsilon_3 u_3 + \epsilon_4 u_4)^4
\, = \, \ \
24 \bigg( \sum_{i=1}^4 u_i^2 \bigg)^2 ,
$$$$
\sum_{i=1}^4 (2u_i)^6 \, + \,
8 \sum_{I} (u_i + \epsilon u_j)^6 \, + \, 
\sum_{\epsilon_i \in \{1,-1\}}
   (u_1 + \epsilon_2 u_2 + \epsilon_3 u_3 + \epsilon_4 u_4)^6
\, = \, \
120 \bigg( \sum_{i=1}^4 u_i^2 \bigg)^3 ,
$$$$
\aligned
6\sum_{i=1}^4 (2u_i)^8 \, + \,
60 \sum_{I} (u_i + \epsilon u_j)^8 \, + \,  
\sum_{II} (2 u_i + \epsilon_j u_j + \epsilon_k u_k)^8  \\
 + \  6 \sum_{\epsilon_i \in \{1,-1\}}
   (u_1 + \epsilon_2 u_2 + \epsilon_3 u_3 + \epsilon_4 u_4)^8
\, &= \, 
5040 \bigg( \sum_{i=1}^4 u_i^2 \bigg)^4 ,
\endaligned
$$$$
\aligned
9\sum_{i=1}^4 (2u_i)^{10} \, + \,
180 \sum_{I} (u_i + \epsilon u_j)^{10} \, + \,  
\sum_{II} (2 u_i + \epsilon_j u_j + \epsilon_k u_k)^{10} \\
  + \  9 \sum_{\epsilon_i \in \{1,-1\}}
   (u_1 + \epsilon_2 u_2 + \epsilon_3 u_3 + \epsilon_4 u_4)^{10}
\, &= \, 
22680 \bigg( \sum_{i=1}^4 u_i^2 \bigg)^5 .
\endaligned
$$
Summations $\sum_{I}$ contain $12$ terms $u_i + \epsilon u_j$,
with $1 \le i < j \le 4$ and $\epsilon = \pm 1$.
Summations $\sum_{II}$ contain $48$ terms of the form
$2u_i + \epsilon_j u_j + \epsilon_k u_k$,
with $i,j,k$ pairwise distinct in $\{1,2,3,4\}$, $j < k$,
and $\epsilon_j, \epsilon_k = \pm 1$.
These four identities provide cubature formulas on $\Bbb S^3$
of sizes and strengths
$$
\aligned
24 \quad &\text{and} \quad 5 \qquad
     \text{(Liouville),} \\
48 \quad &\text{and} \quad 7 \qquad
     \text{(Kempner),} \\
144 \quad &\text{and} \quad 9 \qquad
     \text{(Hurwitz),} \\
144 \quad &\text{and} \quad 11 \qquad
     \text{(Schur).} 
\endaligned
$$
Other examples appear in \cite{Dicks--19} (pages 717--724).

\smallskip

   Observe that Liouville's identity provides 
a design for $\Cal P^{(4)}(\Bbb S^3)$ of size $12$ 
and a spherical $5$-design of size $24$. 
The latter is (up to homothety) a root system of type $D_4$;
in other coordinates, it can also be written as
$$
\sum_{1 \le i < j \le 4, \epsilon \in \{-1,1\}}
(u_i + \epsilon u_j)^4 
\, = \, 
6 \bigg( \sum_{i=1}^4 u_i^2 \bigg)^2 
\tag{L}
$$
(apparently written down first by Lucas in 1876).
This is {\it not} a design for $\Cal P^{(6)}(\Bbb S^3)$,
otherwise there would exist a cubature formula of strength
$7$ on $\Bbb S^3$ of size $24$,
but such cubature formulas cannot have size less than $40$
by Proposition~3.2. 

   Kempner's identity provides a spherical $7$-design in $\Bbb S^3$
of size $48$. 
There is in \cite{HarSl--94} a cubature formula of strength $7$ 
in $\Bbb S^3$ which is of size $46$, conjecturally the optimal size.

   Schur's identity is not optimal 
for cubatures formulas  of strength $11$ on $\Bbb S^3$. 
Indeed, the $120$ vertices of a regular polytope of Schl\"afli symbol
\footnote{
See \cite{Coxet--73}, page 153.
}
$\{3,3,5\}$ provide an antipodal spherical $11$-design
(compare with the lower bound $112$ of Proposition~3.3
for $11$-designs in $\Bbb S^3$). 

\smallskip

   Any orbit in $\Bbb S^3$ for the natural action of the
Coxeter group $H_4$ is a $11$-design,
and there exists a particular orbit of size $1440$
which is indeed a $19$-design.
This construction is apparently due to Salihov (1975). 
See the end of Section 5 in \cite{GoSe--81a},
page 112 of \cite{Rezni--92}, as well as \cite{HarPa--04}.

\bigskip

\noindent {\bf 3.9.\ A digression on Waring problem.} 
Liouville used his identity (see~3.8) and
Lagrange's theorem on the representation of integers as
sums of four squares to show the following claim:

{\it any positive integers is a sum of at most $53$ biquadrates
(= fourth powers).}

\noindent Here is a proof of the claim,
using Lucas' form (L) of Liouville's identity.
Since any positive number is (by Lagrange's theorem) of the form
$6(N_1^2 + N_2^2 + N_3^2 + N_4^2) + r$ with $N_1,\hdots,N_4 \in \Bbb N$
and $r \in \{0,1,2,3,4,5\}$,
it is enough to check that any number of the form $6N^2$ 
is a sum of $12$ biquadrates.
Using Lagrange's theorem again, $N$ can be written as a sum
$n_1^2+n_2^2+n_3^2+n_4^2$ of four squares.
If $n \in \Bbb Z^4$ denotes the vector of coordinates $n_1,n_2,n_3,n_4$,
we have $N = \langle n \mid n \rangle$ and
$$
6N^2 \, = \, 6 \langle n \mid n \rangle ^2 \, = \,
\sum_{1 \le i < j \le 4, \epsilon \in \{-1,1\}}
(n_i + \epsilon n_j)^4  \, = \, 
\text{sum of $12$ fourth powers}
$$
by (L), and this ends the proof.

   Today, we know how to make this bound sharp,
reducing it from $53$ to $19$ \cite{BDD--86a/b},
or even to $16$ if a {\it finite number} of exceptions is allowed
\cite{Daven--39}. More precisely, there are exactly $96$ numbers
which are not sums of $16$ biquadrates, 
and the maximum of them is $13\,792$~; 
their list is shown in \cite{DeHeL--00}.

\bigskip

\noindent {\bf 3.10.\ Lattice cubature formulas.} Here is a variation 
on the construction of lattice designs described in Example~1.16.

Consider a lattice $\Lambda \subset \Bbb R^n$
which is even and unimodular.
For an integer $l \ge 0$ and a harmonic homogeneous polynomial function
$P \in \Cal H^{(2l)}(\Bbb R^n)$, 
the {\bf theta series} is defined by
$$
\Theta_{\Lambda,P} (z) \, = \, 
\sum_{\lambda \in \Lambda} P(\lambda) 
     q^{\frac{1}{2}\langle \lambda \mid \lambda \rangle} \, = \,
P(0) +
\sum_{x \in \Lambda_2} P(x) q + 
\sum_{x \in \Lambda_4} P(x) q^2 + \cdots
$$
where $q = e^{2i\pi z}$,
$z \in \Bbb C$ with $\operatorname{Im}(z) > 0$, and 
$\Lambda_{2j} = 
\{x \in \Lambda \mid \langle x \mid x \rangle = 2j \}$.
It is a modular form for $PSL(2,\Bbb Z)$
of weight $\frac{n}{2} + 2l$,
and a parabolic modular form if moreover $l > 0$.
In particular, 
we have  $\Theta_{\Lambda,P} = 0$ 
whenever $l > 0$ and
$\frac{n}{2} + 2l \le 10$ or $ \frac{n}{2} + 2l = 14$
(see 
% \cite{Hecke--40}, and 
Corollary~3.3 in \cite{Ebeli--94}).

   Suppose now that $n = 8$,
so that $\Lambda$ is a root lattice of type $E_8$.
The previous conside\-rations show that $\Lambda_{2j}$
is a spherical $7$-design 
and a design for $\Cal H ^{(10)}(\Bbb R^8)$
for any $j \ge 1$.
Moreover, for $P \in \Cal H ^{(8)}(\Bbb R^8)$,
the parabolic modular form $\Theta_{\Lambda,P}$ 
is necessarily a constant multiple of
$$
\Delta_{24}(z) \, = \, q^2 \prod_{m=1}^{\infty} (1-q^{2m})^{24}
\, = \, q^2 - 24q^4 + 252q^6 - 1472q^8 + \cdots  
\qquad (q = e^{i\pi z}).
$$
In particular,
$
24 \sum_{x \in \Lambda_2} P(x) + \sum_{x \in \Lambda_4} P(x)
= 0 .
$
Observe that $P(x) = 2^4 P(x/\sqrt 2)$ for $x \in \Lambda_2$
and $P(x) = 2^8 P(x/2)$ for $x \in \Lambda_4$.
It follows that          
$$
\frac{3}{2} \sum_{x \in \Lambda_2} P\left(\frac{x}{\sqrt 2}\right) 
+ \sum_{x \in \Lambda_4} P\left(\frac{x}{2}\right)
\, = \, 0 
$$
for all $P \in \Cal H^{(k)}(\Bbb R^8)$ and $1 \le k \le 11$.
In other words, the first two shells  of $\Lambda$
provide a cubature formula in $\Bbb S^{7}$ of strength $11$,
with underlying set 
$\frac{1}{\sqrt 2} \Lambda_2 \sqcup \frac{1}{2} \Lambda_4$, 
of size
$$
\abs{\Lambda_2} + \abs{\Lambda_4} \, = \, 240 + 2160 \, = \, 2400 ,
$$
and with exactly two different weights,
namely $1/1680$ on $\frac{1}{\sqrt 2} \Lambda_2$
and $1/2520$ on $\frac{1}{2} \Lambda_4$.
%(the ration of the weights being $\frac{3}{2}$).
This size compares favourably with the bounds of Proposition~3.6,
which read here
$$
2 \operatorname{dim}_{\Bbb R}\Cal P^{(5)}(\Bbb R^8) = 1584
\, \le \, N \, \le \,
2 \operatorname{dim}_{\Bbb R}\Cal P^{(10)}(\Bbb R^8)-2 = 38894 .
$$
A similar lower bound shows that 
$\frac{1}{\sqrt 2} \Lambda_2 \sqcup \frac{1}{2} \Lambda_4$
cannot enter a cubature formula of strength $13$ on $\Bbb S^7$.

The previous construction indicates clearly a general procedure.

\bigskip

\noindent {\bf 3.11.\ On Gaussian designs.} Let $l$ be an integer, 
$l \ge 0$, and let $\alpha \in \Bbb R^n$. On the one hand, 
it is classical that
$$
\int_{\Bbb S^{n-1}} \langle \alpha \mid x \rangle ^{2l} d \sigma (x)
\, = \,
\frac{(2l-1)!!}{n(n+2)\cdots (n+2l-2)} \,
\langle \alpha \mid \alpha \rangle^l
$$
(notation: $(2l-1)!! = \prod_{j=1}^l (2l-2j+1)$).
Indeed, the left hand side is clearly a homogeneous polynomial of degree
$2l$ in $\alpha$ which is invariant by the orthogonal group $O(n)$,
and therefore a constant multiple of 
$\langle \alpha \mid \alpha \rangle^l$.
If we apply both sides the Laplacian with respect to $\alpha$,
we find a recurrence relation for the constants, 
and this provides the result.
On the other hand, we have
$$
\frac{1}{\pi^{n/2}} \int_{\Bbb R^n} 
\langle \alpha \mid x \rangle ^{2l} e^{-\norm{x}^2} dx
\, = \, \frac{(2l-1)!!}{2^l} \, 
\langle \alpha \mid \alpha \rangle^l
$$
(calculus, using integration by parts in case $n = 1$).
It follows that
$$
\int_{\Bbb S^{n-1}} \varphi(x) d \sigma (x)
\, = \,
\frac{2^l}{n(n+2)\cdots (n+2l-2)} \ \frac{1}{\pi^{n/2}}
\int_{\Bbb R^n} \varphi (x) e^{-\norm{x}^2} dx
\tag$*$
$$
for all $\varphi \in \Cal P^{(2l)}(\Bbb R^n)$, $l \ge 0$.
[See the footnote in Proposition~1.13.]
Observe that, for $\varphi \in \Cal P^{(2l+1)}(\Bbb R^n)$,
both integrals in $(*)$ vanish for symmetry reasons.

\medskip

   Let now $X$ be a non-empty finite subset of $\Bbb S^{n-1}$ of size $N$.
If $X$ is a spherical $t$-design with $t \in \{1,2,3\}$, the set
$$
X_G \, = \, \sqrt{\frac{n}{2}} X 
$$
is a {\it Gaussian $t$-design.}
Indeed, for $\varphi \in \Cal P^{(2)}(\Bbb R^n)$, we have
$$
\frac{1}{N} \sum_{x \in X_G} \varphi(x) \, = \,
\frac{n}{2N} \sum_{x \in X} \varphi(x) \, = \,
\frac{n}{2} \int_{\Bbb S^{n-1}} \varphi (x) d\sigma(x) \, = \,
\frac{1}{\pi^{n/2}} \int_{\Bbb R^n} \varphi(x) e^{-\norm{x}^2} dx 
$$
(we leave to the reader the verifications with
$\varphi \in \Cal P^{(1)}(\Bbb R^n)$ and
$\varphi \in \Cal P^{(3)}(\Bbb R^n)$).

If  $X$ is now a spherical $5$-design, 
set 
$$
\aligned
\rho_1 \, &= \, \sqrt{\frac{1}{2}(n - \sqrt{2n})} \\
\rho_2 \, &= \, \sqrt{\frac{1}{2}(n + \sqrt{2n})}
\endaligned
\qquad \text{so that} \qquad
\aligned
\frac{1}{2}(\rho_1^2 + \rho_2^2) \, &= \, \hskip.55cm \frac{n}{2} \\
\frac{1}{2}(\rho_1^4 + \rho_2^4) \, &= \, \frac{n(n+2)}{4} .
\endaligned
$$
Then
$$
X_{GG} \, = \, \rho_1 X \sqcup \rho_2 X
\qquad \text{(disjoint union)}
$$
is a {\it Gaussian $5$-design} of size $2N$.
Indeed, for $\varphi \in \Cal P^{(4)}(\Bbb R^n)$, 
we have
$$
\aligned
\frac{1}{2N} \sum_{x \in X_{GG}} \varphi (x) \, &= \,
\frac{\rho_1^{4}}{2N}\sum_{x \in X} \varphi (x) +
\frac{\rho_2^{4}}{2N}\sum_{x \in X} \varphi (x) \, = \,
\frac{1}{2}(\rho_1^4 + \rho_2^4) 
    \int_{\Bbb S^{n-1}} \varphi(x) d\sigma(x) \\
&= \, \frac{1}{2}(\rho_1^4 + \rho_2^4)  \frac{4}{n(n+2)}
     \frac{1}{\pi^{n/2}} \int_{\Bbb R^n} 
     \varphi(x) e^{-\norm{x}^2} dx \\
&= \,    \frac{1}{\pi^{n/2}} \int_{\Bbb R^n} 
     \varphi(x) e^{-\norm{x}^2} dx
\endaligned
$$
(where the last equality is a consequence of the values chosen for
$\rho_1$ and $\rho_2$). Similarly
$$
\frac{1}{2N} \sum_{x \in X_{GG}} \varphi (x) \, = \,
\frac{1}{\pi^{n/2}} \int_{\Bbb R^n} 
     \varphi(x) e^{-\norm{x}^2} dx
$$
for all $\varphi \in \Cal P^{(2)}(\Bbb R^n)$.
We leave to the reader the verifications with $\varphi$
a homogeneous polynomial of odd degree.

\bigskip

We end the present chapter with a proposition and a construction suggesting
one more connection with another subject.

   For an integer $N \ge 1$ and a real number $p \ge 1$,
we denote by $\ell^p(N)$ the classical Banach space of dimension $N$,
with underlying space $\Bbb R^N$ and norm $\norm{x}_p$, where
$$
\norm{\left(x_j\right)_{1 \le j \le N}}_p \, = \,
\root{p} \of{
\smash{
\sum_{1 \le j \le N} \abs{x_j}^p
}\vphantom{\sum}
} .
$$
\vskip.1cm \noindent
A linear mapping from a real vector space to $\ell^p(N)$ 
is said to be {\bf degenerate} 
if its image is in a hyperplane of equation $x_j = 0$ 
for some $j \in \{1,\hdots,N\}$. 
The next proposition is from \cite{LyuVa--93},
and $c_{2l}$ is the constant of our Proposition~1.13.

\bigskip

\proclaim{3.12.\ Proposition} Consider integers
$n \ge 2$, $l \ge 2$, and $N \ge 1$.

   (i) To any cubature formula $(X,W)$ of size $N$ 
for $\Cal P^{(2l)}(\Bbb S^{n-1})$ corresponds 
a non-dege\-nerate isometric embedding
$J_{X,W} : \ell^2(n) \longrightarrow \ell^{2l}(N)$.

   (ii) Conversely, to any non-degenerate isometric embedding 
$J : \ell^2(n) \longrightarrow \ell^{2l}(N)$
corres\-ponds
a cubature formula of size $N$ in $\Cal P^{(2l)}(\Bbb S^{n-1})$.
\endproclaim

\demo{Proof} Let $(X,W)$ be as in $(i)$ and let
$(x^{(1)}, \hdots, x^{(N)})$ be an enumeration of the points in $X$.
The linear mapping
$J_{X,W} : \Bbb R^n \longrightarrow \Bbb R^N$ defined by its coordinates
$$
J_{X,W}(u)_k \, = \, \left( \frac{ W(x^{(k)}) }{c_{2l}} \right)^{1/2l}
   \left\langle x^{(k)} \mid u \right\rangle ,
   \qquad 1 \le k \le N ,
$$
is an isometry $\ell^2(n) \longrightarrow \ell^{2l}(N)$ by 
the argument of Proposition 1.13.
It is non-degenerate since $X^{\perp} = \{0\}$.

   Let $J : \ell^2(n) \longrightarrow \ell^{2l}(N)$ be a non-degenerate
isometry as in $(ii)$. 
For each $k \in \{1,\hdots,N\}$, let $y^{(k)}$ be the unique
vector in $\ell^2(n)$ such that 
$J(u)_k = \left\langle y^{(k)} \mid u \right\rangle$
for all $u \in \ell^2(n)$;
observe that $y^{(k)} \ne 0$, by the non-degeneracy condition on $N$.
Set $x^{(k)} = y^{(k)} / \norm{y^{(k)}}_2$
and $W^{(k)} = c_{2l}\norm{y^{(k)}}_2^{2l}$.
Then $X = \{x^{(1)}, \hdots, x^{(N)}\}$ is a $N$-subset of $\Bbb S^{n-1}$,
the $W^{(k)}$ define a weight $X \longrightarrow \Bbb R_+^*$, 
and $(X,W)$ is a cubature formula for $\Cal P^{(2l)}(\Bbb S^{n-1})$
by Proposition~1.13.
$\square$
\enddemo

\bigskip
\head
{\bf 
4.\ Markov operators
}
\endhead
\medskip

  Given a complex Hilbert space $\Cal H$, 
we denote by $\Cal L (\Cal H)$ the C$^*$-algebra of all bounded linear
operators on $\Cal H$ and by $\Cal U(\Cal H)$ its {\it unitary group,}
consisting of all unitary operators on $\Cal H$.
Given a group $\Gamma$, we denote by $\Bbb C[\Gamma]$ the 
{\it group algebra} of all functions 
$\Gamma \longrightarrow \Bbb C$ of finite support, with multiplication
the convolution; we consider $\Gamma$ as a subset of $\Bbb C[\Gamma]$,
by identifying $\gamma \in \Gamma$ with the function of value $1$ on
$\gamma$ and $0$ elsewehere. 
A {\it unitary representation} of a group $\Gamma$ in a Hilbert space 
$\Cal H$ is a group homomorphism
$\pi : \Gamma \longrightarrow \Cal U(\Cal H)$;
it extends to a morphism of algebra 
$\Bbb C[\Gamma] \longrightarrow \Cal L (\Cal H)$, 
denoted by $\pi$ again, and defined by 
$\pi(f) = \sum_{\gamma \in \Gamma} f(\gamma) \pi(\gamma)$.
Recall that the norm of a unitary operator is $1$, so that
$\norm{\pi(f)} \le \sum_{\gamma \in \Gamma} \abs{f(\gamma)}$.
Recall also that, if $f(\gamma^{-1}) =  \overline{f(\gamma)}$
for all $\gamma \in \Gamma$, the operator $\pi(f)$ is self-adjoint;
in particular, the spectrum of $\pi(f)$ is then a closed subset of the
real interval $[-\norm{\pi(f)} , \norm{\pi(f)}]$.

   Up to minor terminological changes, this holds for
a {\it real} Hilbert space $\Cal H$, 
its {\it ortho\-gonal group} $\Cal O (\Cal H)$, 
{\it orthogonal representations} $\Gamma \longrightarrow \Cal O (\Cal H)$,
and corresponding morphisms 
$\Bbb R[\Gamma] \longrightarrow \Cal L (\Cal H)$.

\bigskip

\proclaim{4.1.\ Definition} Let $\Gamma$ be a group, 
$S$ a finite subset of $\Gamma$, and $W : S \longrightarrow \Bbb R_+^*$
a weight function; 
set $M_{S,W} = \sum_{s \in S}W(s) s \in \Bbb R [\Gamma]$.
The {\bf Markov operator} associated to the pair~$(S,W)$ 
and an orthogonal (or unitary) representation
$\pi : \Gamma \longrightarrow \Cal U(\Cal H)$ is the operator
$$
\pi(M_{S,W}) \, = \,  \sum_{s \in S} W(s) \pi (s)
$$
on $\Cal H$.
\endproclaim

   Markov operators (often with constant weights) play an important role
in the study of random walks on Cayley graphs of groups \cite{Keste--59},
and more generally in connection with 
unitary representations of groups. See for example 
\cite{BarGr--00},
\cite{BeVaZ--97},
\cite{GaJaS--99},  
\cite{GriZu--01},
\cite{GriZu--02},
\cite{GriZu},
\cite{HaRV1--93}, and
\cite{HaRV2--93}.

\bigskip

  Consider a sigma-finite measure space $(\Omega,\sigma)$ and a group $G$
acting on $\Omega$ by measurable transformations preserving the measure
class of $\sigma$. The corresponding unitary representation $\rho$ of $G$
on $L^2(\Omega,\sigma)$ is defined by
$$
(\rho(g)\varphi)(\omega) \, = \, 
\sqrt{\frac{d g \sigma}{d \sigma}} \ \varphi(g^{-1}\omega)
$$
for $g \in G$, $\varphi \in L^2(\Omega,\sigma)$, 
and $\omega \in \Omega$, where
$\frac{d g \sigma}{d \sigma}$ denotes the appropriate Radon-Nikodym
derivative. 
Given a mapping $s : \Omega \longrightarrow G$, 
$\omega \longmapsto s_{\omega}$,
a finite subset $X$ of $\Omega$,
and a weight function $W : X \longrightarrow \Bbb R_+^*$, 
we have an element
$M_{X,W} = \sum_{x \in X}W(x)s_x$
in the real group algebra of the subgroup of $G$
generated by $s(X)$.
Given moreover a representation $\pi$ of $G$ 
in a Hilbert space $\Cal H_{\pi}$, 
we have a Markov operator
$$
\pi(M_{X,W}) \, = \, \sum_{x \in X} W(x) \pi(s_x) \, \in \,
\Cal L (\Cal H_{\pi}) .
$$

   The particular case studied in \cite{Pache--04} is that of the
orthogonal group $G = O(n)$ acting on the unit sphere $\Omega = \Bbb
S^{n-1}$,    with $\Cal H_{\pi}$ one of the finite-dimensional spaces
introduced in Example~1.11 above. 
In this case, the image of a point $x \in \Bbb S^{n-1}$
by the mapping $s$ is the orthogonal reflection
of $\Bbb R^n$ which fixes the hyperplane orthogonal to $x$.

\bigskip

\proclaim{4.2.\ Proposition} Let $n \ge 2$, $l \ge 0$ be integers
and let $\pi^{(l)}$ denote the natural representation of the group $O(n)$
in the space $\Cal H ^{(l)}(\Bbb S^{n-1})$
of harmonic homogeneous polynomials of degree $l$.
For each $x \in \Bbb S^{n-1}$, let $s_x \in O(n)$ denote the reflection
of $\Bbb R^n$ that fixes the hyperplane $x^{\perp}$ of $\Bbb R^n$.

   If $(X,W)$ is a cubature formula of strength $2l$ on $\Bbb S^{n-1}$,
then
$$
\pi^{(l)}(M_{X,W}) \, = \,
\frac{n-2}{2l+n-2} \operatorname{id}^{(l)}
$$
where $\operatorname{id}^{(l)}$ denotes the identity operator on
$\Cal H ^{(l)}(\Bbb S^{n-1})$.
\endproclaim

\demo\nofrills{Proof, repeated from \cite{Pache--04} \usualspace}
(see also \cite{HarVe--01}).
For $x \in \Bbb R^n$, we define the operator $\tilde{s}_x$ on 
$\Bbb R^n$ by 
$\tilde{s}_x(u) = \langle x \mid x \rangle u
- 2 \langle x \mid u \rangle x$.
Note that $\tilde{s}_x$ is selfadjoint, 
that $\tilde{s}_{\lambda x} = \lambda^2\tilde{s}_x$ 
for $\lambda \in \Bbb R$,
and that $\tilde{s}_x = s_x$ if $x \in \Bbb S^{n-1}$.\par

   For $k \ge 0$,  define the selfadjoint operator
$$
\overline{M}^{(k)} \, = \,
\int_{\Bbb S^{n-1}} \pi^{(k)}(s_{\omega}) d\sigma(\omega)
$$
on $\Cal H^{(k)}(\Bbb R^n)$. For $g \in O(n)$, we have
$s_{g(\omega)} = gs_{\omega}g^{-1}$. 
Since the measure $\sigma$ is $O(n)$-invariant 
and the representation $\pi^{(k)}$ is
irreducible, it follows from Schur's lemma that
$\overline{M}^{(k)}$ is a homothety. By a simple computation, we
obtain the trace of $\pi^{(k)}(s_{\omega})$, which is independent
of $\omega$, and therefore the trace of $\overline{M}^{(k)}$. It
follows that
$$
\overline{M}^{(k)} \, = \, \frac{n-2}{2k+n-2} \operatorname{id}^{(k)} 
\tag1
$$
for all $k \ge 0$. \par

   For $\varphi_1,\varphi_2 \in \Cal H^{(k)}(\Bbb R^n)$
and $u \in \Bbb R^n$, set
$$
\Psi^{(k)}(\varphi_1,\varphi_2)(u) 
\, = \,
\langle \varphi_1 \mid \varphi_2 \circ \tilde{s}_u \rangle
\, = \, 
\int_{\Bbb S^{n-1}} \varphi_1(\omega)\varphi_2(\tilde{s}_u(\omega))
d\sigma(\omega) .
$$
It is easy to check that the function
$\Psi^{(k)}(\varphi_1,\varphi_2)$ is polynomial, homogeneous of
degree $2k$, and depends symmetrically on $\varphi_1,\varphi_2$.
Thus, we have a linear mapping
$$
\Psi^{(k)} \, : \, 
\operatorname{Sym}^2(\Cal H^{(k)}(\Bbb R^n))
\ \longrightarrow \ 
\Cal P^{(2k)}(\Bbb R^n) .
$$
Moreover
$\int_{\Bbb S^{n-1}} \Psi^{(k)}(\varphi_1,\varphi_2)(\omega) d\sigma(\omega)
\, = \, 
\left\langle \varphi_1 \Bigm| \overline{M}^{(k)} (\varphi_2) \right\rangle$.
\par

   Consider now a non-empty finite subset $X$ of $\Bbb S^{n-1}$,
a weight $W : X \longrightarrow \Bbb R_+^*$, 
and an integer $l \ge 0$. Then
$$
\sum_{x \in X} W(x) \varphi (x) \, = \,
\int_{\Bbb S^{n-1}} \varphi (\omega) d\sigma(\omega)     
\qquad \text{for all $\varphi$ in the image of $\Psi^{(l)}$}
\tag2
$$
if and only if
$$
\sum_{x \in X}
W(x) \langle \varphi_1 \mid \varphi_2 \circ s_x \rangle
\, = \,
\left\langle \varphi_1 \Bigm| \overline{M}^{(l)} (\varphi_2) \right\rangle
\qquad
\text{for all $\varphi_1,\varphi_2 \in \Cal H^{(l)}(\Bbb R^n)$,}
$$
namely, 
by the equality $\pi^{(l)}(M_{X,W})\varphi_2 = 
\sum_{x \in X}W(x) \varphi_2 \circ s_x$ defining $\pi^{(l)}(M_{X,W})$
and by~$(1)$,  
if and only if
$$
\pi^{(l)}(M_{X,W}) \, = \, \overline{M}^{(l)} \, = \,
\frac{n-2}{2l+n-2} \operatorname{id}^{(l)}.    
\tag3
$$

   If $(X,W)$ is a cubature formula of strength $2l$, then
$(2)$ holds for all $\varphi \in\Cal P^{(2l)}(\Bbb R^n)$,
and a fortiori for all $\varphi$ in the image of $\Psi^{(l)}$, so
that the proposition is proved. $\square$
\enddemo

\bigskip

   It is remarkable that a converse holds for $n > 2$.

\proclaim{4.3.\ Proposition} Let $n \ge 3$, $l \ge 0$ be integers
and let $\pi^{(l)}$ be as in the previous proposition.
Let $X$ be a non-empty antipodal finite subset of $\Bbb S^{n-1}$
and let $W : X \longrightarrow \Bbb R_+^*$ be a symmetric weight.

   If $\pi^{(l)}(M_{X,W})$ is a homothety, 
then $(X,W)$ is a
cubature formula of strength 
$2l+1$.
\endproclaim

\demo{Proof} See \cite{Pache--04}. $\square$
\enddemo

\bigskip

More generally, consider a Riemannian symmetric pair $(G,H)$ 
where $G$ is a compact Lie group acting on $\Omega = G/H$, 
and the $G$-invariant probability measure $\sigma$ on $\Omega$.
Let $s = s^{\operatorname{symm}} : \Omega \longrightarrow G$ 
be the mapping which associates to a point $x \in G/H$ 
the symmetry of $\Omega$ fixing $x$. 
(For spheres, observe that $s^{\operatorname{symm}}(x)$
is {\it minus} the reflection $s_x$ fixing $x^{\perp}$.) 
There is an orthogonal decomposition
$$
L^2(\Omega,\sigma) \, = \, \bigoplus_{\lambda \in \Lambda} V^{\lambda}
$$
in irreducible $G$-spaces.
For a finite subset $X$ of $\Omega$ and a weight function $W$ on $X$,
there is an analogue of Proposition~4.2  concerning
spaces $V^{\lambda}$ for which $(X,W)$ is a cubature formula
and for which the Markov operator
$\pi^{(\lambda)}(M_{X,W}) \in \operatorname{End}(V^{(\lambda)})$ 
is a constant multiple of the identity
\cite{Pache, in preparation}.
For designs in $G/H$ a Grassmannian, 
see \cite{BaCoN--02} and \cite{BaBaC--04}.

\bigskip

   Consider a non-empty finite subset $S$ of $O(n)$ which is symmetric
($S^{-1} = S$) and denote by $\Gamma_S$ the subgroup of $O(n)$
it generates.
\par
\centerline{
\bf{From now on, we assume that
$W(s) = \abs{S}^{-1}$ for all $s \in S$} 
}
\par
\centerline{
\bf{and we set $M_S = \abs{S}^{-1}\sum_{s \in S} s \in \Bbb R[\Gamma_S]$.}
}
\par\noindent
Propositions~4.2 and 4.3 relate the spectra of Markov operators
$\pi^{(k)}(M_S)$ for some small values of $k$ to the cubature
properties of $(X,W)$. 
It happens that the spectra for $k \to \infty$
of the operators $\pi^{(k)}(M_S)$
are also important, as shown by Lubotzky, Phillips, and Sarnak
\cite{LuPhS--86, 87}
(see also \cite{ColdV--88}). 
We formulate the weak Observation~4.4 before the stronger Proposition~4.5.

Let $\pi_0$ denote the natural representation of $O(n)$ in the space
$L^2_0(\Bbb S^{n-1},\sigma)$ of $L^2$-functions of zero average. 
The orthogonal sum
$L^2_0(\Bbb S^{n-1},\sigma) = 
\bigoplus_{k = 1}^{\infty} \Cal H^{(k)}(\Bbb S^{n-1})$
provides an orthogonal decomposition
$
\pi_0 \, = \, \bigoplus_{k=1}^{\infty} \pi^{(k)}
$
into irreducible subrepresentations (see, e.g., \cite{SteWe--71}), so that
$$
\norm{\pi_0(M_S)} \, = \, \sup_{k \ge 1} \norm{\pi^{(k)}(M_S)} .
$$

\bigskip

\proclaim{4.4.\ Observation} For any finite symmetric subset $S$ of $O(n)$,
we have
$$
1 \, \ge \, \norm{\pi_0(M_S)} \, \ge \, \frac{1}{\sqrt{\abs{S}}} .
$$
\endproclaim

\demo{Proof} (Compare with Proposition~2.3.2 in \cite{Sarna--90}.)
For all $s \in S$, we have $\norm{\pi_0(s)} = 1$, since the
representation $\pi_0$ is orthogonal. The upper bound follows.

Choose $\omega \in \Bbb S^{n-1}$ such that $s(\omega) \ne t(\omega)$
for all $s,t \in S \cup \{id\}$, $s \ne t$.
There exists a neighbourhood $U$ of $\omega$ in $\Bbb S^{n-1}$
such that $s(U) \cap t(U) = \emptyset$ 
for all $s,t \in S \cup \{id\}$, $s \ne t$,
and a function $\varphi \in L^2_0(\Bbb S^{n-1},\sigma)$
of norm $1$ supported in $U$.
Since $\norm{\pi_0(M_S)\varphi}^2 = \abs{S}^{-1}$,
we have $\norm{\pi_0(M_S)} \ge \abs{S}^{-1/2}$.
$\square$
\enddemo

\bigskip

\proclaim{4.5.\ Proposition} Let $S$ be a symmetric finite  subset of
$O(3)$. Then
$$
\norm{\pi_0(M_S)} 
\, \ge \, \limsup_{k \to \infty} \norm{\pi^{(k)}(M_S)}
\, \ge \, \frac{2\sqrt{\abs{S}-1}}{\abs{S}} .
$$
\endproclaim

 \medskip

  The proofs of \cite{LuPhS--86} and \cite{ColdV--88} are
written up for the case of a subset $S$ of $SO(n)$. 
The present generalization to $O(n)$ is rather straightforward. 
We isolate part of the proof in the following lemma

\bigskip

\proclaim{4.6.\ Lemma} (i)
Let $g \in SO(3)$ be a rotation of angle $\theta_g \in [0,\pi]$.
Then
$$
\operatorname{trace}\left( \pi^{(k)}(g) \right) \, = \, \left\{
\aligned
\frac{\sin [(2k+1)\theta_g/2]}{\sin[\theta_g/2]}
   \quad &\text{if} \quad \theta_g \ne 0 \\
2k+1
   \hskip1cm &\text{if} \quad \theta_g = 0 .
\endaligned
\right.
$$

  (ii) Let $g \in O(3)$ have eigenvalues $\exp (\pm i \theta_g)$ and $-1$.
Then
$$
\operatorname{trace}\left( \pi^{(k)}(g) \right) \, = \, \left\{
\aligned
\frac{\cos [(2k+1)\theta_g/2]}{\cos[\theta_g/2]}
   \quad &\text{if} \quad \theta_g \ne \pi \\
(-1)^k
   \hskip1cm &\text{if} \quad \theta_g = \pi .
\endaligned
\right.
$$
\endproclaim

\demo{Proof} $(i)$ The eigenvalues of $g \in SO(3)$ are 
$\exp(\pm i \theta_g)$ and $1$. 
The space $\Cal P^{(k)}(\Bbb S^ 2)$ has a linear basis
of eigenvectors of the transformation induced by $g$
of the form
$$
\left\{ x^a y^b z^{k-a-b} \mid a,b \ge 0 \ \text{and} \ a+b \le
k\right\} .
$$
For $j \in \{0,\hdots,k\}$, the trace of the linear
endomorphism defined by $g$ on the linear span of 
$\left\{ x^a y^b z^{k-a-b}  \mid a,b \ge 0 \ \text{and} \ a+b = j\right\}$
is
$$
\sigma_j(g) \, = \, e^{ij\theta_g} + e^{i(j-2)\theta_g} + \cdots +
      e^{-ij\theta_g} 
     \, = \, 
    \frac{ \sin[(j+1)\theta_g]}{\sin [\theta_g]}
$$
so that the trace of the linear endomorphism defined by $g$ on
$\Cal P^{(k)}(\Bbb S^2)$ is $\sum_{j=0}^k \sigma_j (\theta_g)$.
Since $\Cal P^{(k)}(\Bbb S^2) = 
\bigoplus_{j \ge 0, k-2j \ge 0} \Cal H^{(k-2j)}(\Bbb S^2)$, we have
$\operatorname{trace}\left( \pi^{(k)}(g) \right) \, = \,
\sigma_k (\theta_g) + \sigma_{k-1}(\theta_g)$
and the formula of $(i)$ follows.

   $(ii)$ Similarly, for $g \in O(3)$, $g \notin SO(3)$, the trace of the
linear endomorphism defined by $g$ on $\Cal P^{(k)}(\Bbb S^2)$ is
$\sum_{j=0}^k (-1)^{k-j} \sigma_j (\theta_g )$ so that
$$
\operatorname{trace}\left( \pi^{(k)}(g) \right) \, = \, 
\sigma_k(\theta_g) - \sigma_{k-1}(\theta_g) .
$$
(The {\it only} part of the proof which is not explicitely in
\cite{LuPhS--86} is $(ii)$.)
$\square$
\enddemo

\medskip

\demo{Proof of Proposition~4.5}
{\it Step one.}
Recall that $\Gamma_S$ is the subgroup of~$O(3)$ generated by~$S$. 
For each integer $N \ge 0$, let $W_N$ denote the number of words
\footnote{
Words which need not be reduced in any sense.
}
in letters of $S$, of length $N$.
To each word $w \in W_N$ corresponds naturally an orthogonal transformation
$g = g(w) \in \Gamma_S \subset O(3)$.
We have $(M_S)^N = \abs{S}^{-N} \sum_{w \in W_N} g(w)$
and therefore
$$
\abs{S}^N \operatorname{trace}\left( \pi^{(k)}((M_S)^N) \right) \, = \, 
\sum_{w \in W_N \atop g(w) \in SO(3)} 
   \frac{\sin [(2k+1)\theta_g/2]}{\sin[\theta_g/2]} \, + \, 
\sum_{w \in W_N \atop g(w) \notin SO(3)} 
   \frac{\cos [(2k+1)\theta_g/2]}{\cos[\theta_g/2]}
$$
by the previous lemma.

  Let $m_N$ denote the quotient by $\abs{S}^N$ of the number of words 
$w \in W_N$ such that $g(w) = 1 \in O(3)$. We have 
$$
\lim_{k \to \infty} \frac{1}{2k+1}
\operatorname{trace}\left( \pi^{(k)}((M_S)^N) \right) \, = \, m_N
$$
for each $N \ge 0$.
(This is Theorem~1.1 in \cite{LuPhS--86}; observe however the change in
notation: $m_N$ there is $\abs{S}^N$ times what $m_N$ is here and in
\cite{Keste--59}.)

\medskip

{\it Step two.} Consider the Cayley graph
\footnote{
Two elements $\gamma_1,\gamma_2 \in \Gamma_S$
viewed as vertices of $\operatorname{Cay}_S$
are joined by an edge if $\gamma_1^{-1}\gamma_2 \in S$.
We assume that $1 \notin S$, so that the Cayley graph
is simple (= without loops), of degree $\abs{S}$.
}
$\operatorname{Cay}_S$ of
$\Gamma_S$ with respect to $S$, the left regular representation
$\lambda_S$ of $\Gamma_S$, the corresponding Markov operator 
$\lambda_S(M_S)$, and its spectral measure $\mu_{S}$.
Kesten has shown that \par

  $m_N = \int_{\Bbb R} t^N d \mu_{S}(t)$ for all $N \ge 0$, \par

  $\norm{\lambda_S(M_S)} = \limsup_{N \to \infty} 
      \root  {N} \of{m_N}$, \par

   $\norm{\lambda_S(M_S)} \ge 2 \sqrt{ \abs{S} - 1} / \abs{S}$
   with equality if and only if $\operatorname{Cay}_S$ is a tree.

\par\noindent
See \cite{Keste--59}.

\medskip

{\it Step three.} For each $k \ge 0$, let $\mu_{S,k}$ denote the
quotient by $\operatorname{dim}_{\Bbb R}\Cal H^{(k)}(\Bbb S^2) = 2k+1$
of the spectral measure of $\pi^{(k)}(M_S)$;
this can be written
$$
\mu_{S,k} \, = \, \frac{1}{2k+1} \sum_{j=1}^{2k+1} \delta(\lambda_{k,j})
$$
where $\delta(\lambda)$ denotes a Dirac measure of support $\lambda$
and where $\left( \lambda_{k,j} \right)_{1 \le j \le 2k+1}$ 
are the eigenvalues of the endomorphism $\pi^{(k)}(M_S)$ of
$\Cal H^{(k)}(\Bbb S^2)$. Let $N \ge 0$; from the definition of
$\mu_{S,k}$, we have
$$
\frac{1}{2k+1} \operatorname{trace}\left( \pi^{(k)}((M_S)^N) \right)
\, = \,
\int_{\Bbb R} t^N d \mu_{S,k} \, ;
$$
from the two previous steps, we have
$$
\lim_{k \to \infty} \frac{1}{2k+1}
   \operatorname{trace}\left( \pi^{(k)}((M_S)^N) \right)
\, = \, m_N \, = \,
\int_{\Bbb R} t^N d \mu_{S} \, ;
$$
it follows that the sequence of measures $\left(\mu_{S,k}\right)_{k \ge 1}$
(all of mass at most $1$ since
$\norm{\pi^{(k)}(M_S)} \le 1$ and $\norm{\lambda_S(M_S)} \le 1$)
converges vaguely to $\mu_S$, and in particular that
$$
\limsup_{k \to \infty} \norm{\pi^{(k)}(M_S)} \, \ge \,
\frac{2 \sqrt{ \abs{S}-1}}{\abs{S}} \ .
$$
$\square$ 
\enddemo

\bigskip

\proclaim{4.7.\ Theorem (Lubotzky-Phillips-Sarnak)} 
For each prime $p$, there exists a subset $S_0$ of $SO(3)$ of $p+1$
elements such that, if $S = S_0 \cup (S_0)^{-1}$, 
$$
\norm{\pi^{(k)}(M_S)} \, \le \,
\frac{2 \sqrt{ \abs{S}-1}}{\abs{S}} 
\, = \, \frac{\sqrt{2p+1}}{p+1}
$$
for all $k \ge 1$.
\endproclaim

\bigskip

\noindent
{\bf 4.8.\ Open problem.}
It is a natural problem, for each $l \ge 0$, to look for a finite set 
$S \subset O(3)$ of reflections 
(depending on $l$) such that the  following conditions are
fulfilled:
$$
\pi^{(k)}(M_S) \, = \, \frac{1}{2k+1} \operatorname{id}^{(k)}
\qquad \text{for} \quad k \in \{1,\hdots,l\}
\tag*
$$
(see Propositions~4.2 and 4.3),
$$
\norm{\pi^{(k)}(M_S)} \, \le \, 
\frac{2 \sqrt{ \abs{S}-1}}{\abs{S}}
\qquad \text{for $k$ large enough} 
\tag**
$$
(see Proposition~4.5 and Theorem~4.7),
and $\abs{S}$ as small as possible.

Observe that (*) is a way to write that all eigenvalues of
$\pi^{(k)}(M_S)$ are equal to $1/(2k+1)$, for $k \in \{1,\hdots,l\}$,    
and that (**) is a  bound on the eigenvalues of $\pi^{(k)}(M_S)$,
for large values of~$k$.

\bigskip
\head
{\bf 
5.\ Reflection groups
}
\endhead
\medskip

   Let $V$ be a Euclidean space. 
The {\it reflection} $s_x \in O(V)$ associated to $x \in V$, 
$x \ne 0$, is defined as above by
$$
s_x(y) \, = \, y - 
2 \frac{\langle x \mid y \rangle}{\langle x \mid x \rangle}x
\quad \text{for all} \quad y \in V .
$$
Observe that $s_{x'} = s_x$ if and only if $\Bbb R x' = \Bbb R x$.
In this chapter, a {\bf reflection group} is a subgroup of $O(V)$
generated by a finite set of reflections; it {\it need not} act
properly on~$V$, 
contrarily to what is assumed most often in \cite{Bourb--68}. 
Thus, any finite subset $X$ of $V \setminus \{0\}$
defines a reflection group 
$\Gamma_X  \subset O(V)$ generated by $\left(s_x\right)_{x \in X}$. \par

  In particular, let $\Lambda$ be a lattice in $V$. 
We assume that $\Lambda$ is integral, 
so that $\Lambda$ is the disjoint union of the origin 
and of its non-empty  shells
$\Lambda_m, m \ge 1$, as defined in Item~1.15.
Let $\Gamma_{\Lambda,m}$ denote the reflection group generated by
$\{s_x \mid x \in \Lambda_m \}$.
For example, if $\Lambda$ is a root lattice 
(namely an even  lattice $\Lambda$ generated by $\Lambda_2$),
then $\Gamma_{\Lambda,2}$ is a finite group,
and indeed a direct product of Coxeter groups of types
$A_n$ ($n \ge 1$), $D_n$ ($n \ge 4$), and $E_n$ ($n = 6,7,8$);
see for example \cite{Ebeli--94}.
In most cases however, $\Gamma_{\Lambda,m}$ is an infinite group.

   Let us specialize to $m$ a power of $2$.
The scalar product on $V$ defines 
a symmetric $\Bbb Z$-bilinear form on $\Lambda$ 
which extends to a symmetric $\Bbb Z[1/2]$-bilinear form
$$
\beta \, : \, 
(\Lambda \otimes \Bbb Z[1/2]) \times (\Lambda \otimes \Bbb Z[1/2])
\longrightarrow \Bbb Z[1/2].
$$ 
We denote by $O(\Lambda,\Bbb Z[1/2])$ the orthogonal group of $\beta$. 
Observe that $\Gamma_{\Lambda,m} \subset O(\Lambda,\Bbb Z[1/2])$.

   In case of the cubical lattice $\Bbb Z^n \subset \Bbb R^n$,
we write $O(n,\Bbb Z[1/2])$ rather than $O(\Bbb Z^n,\Bbb Z[1/2])$.
Seve\-ral lattices $\Lambda$ 
define the same group $O(n,\Bbb Z[1/2])$; 
a sufficient condition for the group
$O(\Lambda,\Bbb Z[1/2])$ to be isomorhic to $O(n,\Bbb Z[1/2])$
is that $2^{k}\Bbb Z^n \subset \Lambda  \subset 2^{-l}\Bbb Z^n$
for some $k,l \ge 0$.
In particular, 
if $\Lambda \subset \Bbb R^8$ is a root lattice of type $E_8$,
then $O(\Lambda,\Bbb Z[1/2]) = O(8,\Bbb Z[1/2])$,
since $\Lambda$ and $\Bbb Z^8$ have a common sublattice of index $2$
(which is a root lattice of type $D_8$).
Here are some properties of the group $O(n,\Bbb Z[1/2])$:
\roster
\item"(o)"
it is generated by reflections for any $n$;
\item"(i)"
it is infinite if and only if $n \ge 5$;
\item"(ii)" 
it is naturally 
\footnote{
More generally, 
if $\Bbb P$ is a finite set of distinct primes 
and if $N$ denotes their product,
$O(n,\Bbb Z[1/N])$ is naturally a discrete cocompact subgroup of
the product of $O(n)$ with $\prod_{p \in \Bbb P} O(n,\Bbb Q_{p})$,
and therefore also a discrete cocompact subgroup 
of $\prod_{p \in \Bbb P} O(n,\Bbb Q_{p})$. 
It follows that $O(n,\Bbb Z[1/N])$ is finitely presented 
\cite{Serre--71}. 
}
a discrete cocompact subgroups of the $2$-adic group
$O(n,\Bbb Q_2)$;
\item"(iii)"
it is finitely presented;
\item"(iv)"
it is virtually torsion free, and more precisely
$\operatorname{Ker}\big(O(n,\Bbb Z[1/2]) \longrightarrow 
O(n,\Bbb Z/3)\big)$ is torsion free;
\item"(v)"
it has finite-type homology---more precisely, 
for any Noetherian ring $A$ and any $i \ge 0$,  
the $A$-module $H_i(O(n,\Bbb Z[1/2]),A)$ is finitely generated; 
\item"(vi)"
it is  of virtual cohomological dimension
$\frac{1}{2}(n-\inf_{m \in \Bbb Z}\abs{n-8m})$, 
and the group $H_i(O(n,\Bbb Z[1/2]),\Bbb Q)$ 
is zero in dimensions neither $0$ and nor maximal;
\item"(vii)"
the natural inclusion of $O(n,\Bbb Z[1/2])$ into the compact orthogonal
group $O(n)$ has dense image if and only if the group $O(n,\Bbb Z[1/2])$ is
infinite.
\endroster
Several of these properties are straightforward consequences of (ii).
Property (vii) follows from a strong approximation theorem
due to Kneser \cite{Knese--66}.
For much more on these groups, see \cite{Colli}.

\bigskip

The next proposition shows a remarkable generating set of 
$O(8,\Bbb Z[1/2])$; it is an unpublished result of B.~Venkov.

\bigskip

\proclaim{5.1.\ Proposition (B.\ Venkov)} 
Let $L$ be a root lattice of type $E_8$. 
Then the group $O(8,\Bbb Z[1/2])$
is generated by $\Gamma_{L,4}$ and a reflection with respect to one root,
and also by the finite Weyl group $\Gamma_{L,2}$ of type $E_8$ and a
reflection with respect to one element of $L_4$.
\endproclaim

\demo{Proof} Let $\Cal L$ be the set of even unimodular lattices $M$
in $V$ which are  such that $M \cap L$
is of index a power of $2$ in both $L$ and $M$
(the same power since $L$ and $M$ are both unimodular).
Then $\Cal L$ is a metric space for the distance $\delta$ defined by
$$
\delta(M,N) \, = \,
\log_2([M : M\cap N]) \, = \, \log_2([N : M\cap N]) .
$$
It is a theorem of Kneser \cite{Knese--57, page 242}
that $\Cal L$ is connected by steps of length one;
more precisely, given $M,N \in \Cal L$ with $d = \delta(M,N)$,
there exists a sequence 
$M=M^{(0)},M^{(1)},\hdots,M^{(d)}=N$ in $\Cal L$
such that $\delta(M^{(j-1)},M^{(j)}) = 1$ 
for all $j \in \{1,\hdots,d\}$.

Set $G = O(8,\Bbb Z[1/2])$;
this group operates naturally on $\Cal L$.
Here is the synopsis of the proof: for $g \in G$,
set $M = g(L) \in \Cal L$ and $d = \delta(L,M)$;
we show by induction on $d$ 
that $g$ is a product of appropriate reflections. 
More precisely, 
for $M^{(0)}=L, M^{(1)}, \hdots, M^{(d)} = M$ as above, 
we check that $M^{(j)} = s_x(M^{(j-1)})$
for some $x \in M^{(j-1)}$ with $\langle x \mid x \rangle = 4$.
\enddemo

\demo{Step one:
$\langle \Gamma_{L,2},\Gamma_{L,4}\rangle = 
\langle s_r,\Gamma_{L,4}\rangle =
\langle \Gamma_{L,2},s_x\rangle$
for any $r \in L_2$ and $x \in L_4$}

\par\noindent
(For a subset $S \subset O(V)$, we denote by
$\langle S \rangle$ the subgroup of $O(V)$ generated by $S$.) \par

   It is known that the group of all automorphisms of $L$
coincides with the group $\Gamma_{L,2}$ generated
by the symmetries with respect to the roots.
For $g \in \operatorname{Aut}(L)$ and $x \in L$, $x \ne 0$,
we have
$$
s_{g(x)} \, = \, gs_x g^{-1} . \tag$*$
$$
As $\operatorname{Aut}(L)$ acts transitively on $L_4$, we have
$\langle \Gamma_{L,2},\Gamma_{L,4}\rangle = 
\langle \Gamma_{L,2},s_x\rangle$
for any  $x \in L_4$.
(It is known that, more generally, 
$\operatorname{Aut}(L)$ acts transitively on each of
$L_2$, $L_4$, $L_6$, the complement of $2L_2$ in $L_8$,
$L_{10}$, and $L_{12}$; see page 122 of \cite{ConSl--99}.)

   Relations $(*)$ show that $\Gamma_{L,2} \cap \Gamma_{L,4}$
is a normal subgroup of $\Gamma_{L,2}$.
Let $e_1,\hdots,e_8$ be an orthonormal basis of $V$ such that
$r=e_1+e_2$ and $r'=e_1-e_2$ are in $L_2$ 
(see page 268 of \cite{Bourb--68});
then $x=2e_1$ and $x' = 2e_2$ are in $L_4$.
A straightforward computation shows that
$\operatorname{id}_V \ne s_rs_{r'} = s_xs_{x'} \in
\Gamma_{L,2} \cap \Gamma_{L,4}$.
As $\Gamma_{L,2}$ is almost simple
(by Exercise 2 of \S \ VI.4 in \cite{Bourb--68}), 
it follows that $\Gamma_{L,2} \cap \Gamma_{L,4}$ 
is a subgroup of index
\footnote{
We do not know if $\Gamma_{L,2} \cap \Gamma_{L,4}$
is the whole of $\Gamma_{L,2}$ or if it is the
subgroup of elements which are products of
even numbers of reflections.
}
at most $2$ in $\Gamma_{L,2}$, so that $\Gamma_{L,2}$ is generated by
$\Gamma_{L,2} \cap \Gamma_{L,4}$ and $s_r$ for any $r \in L_2$.
A fortiori $\langle \Gamma_{L,2},\Gamma_{L,4}\rangle = 
\langle s_r, \Gamma_{L,4} \rangle$
for any  $r \in L_2$.
\enddemo

\demo{Step two: a reminder on neighbours}

   For this step, $V$ can be a Euclidean space of any dimension
$n \equiv 0 \pmod 8$ and $M$ any even unimodular lattice in $V$.
We view $\overline M = M/2M$ as a vector space of dimension $n$
over the prime field $\Bbb F_2$.
There is a nondegenerate quadratic form
$q : \overline M \longrightarrow \Bbb F_2$
defined by
$q(\overline z) = \frac{1}{2}\langle z \mid z \rangle \pmod{2}$
for $z \in M$ representing
$\overline z \in \overline M$.
For each $z \in M$ such that $z \notin 2M$ and
$\langle z \mid z \rangle \equiv 0 \pmod 4$, set
$$
M_z \, = \, \left\{ m \in M \mid 
\langle m \mid z \rangle \equiv 0 \pmod{2} \right\}
\quad \text{and} \quad
M^z \, = \, M_z \sqcup (\frac{1}{2}z + M_z) .
$$
Then $M^z$ is an integral unimodular lattice in $V$,
and $M \cap M^z$ is of index $2$ in both $M$ and $M^z$.
Moreover:
\roster
\item"(i)" $M^z$ is even if 
$\langle z \mid z \rangle \equiv 0 \pmod 8$
and odd if
$\langle z \mid z \rangle \equiv 4 \pmod 8$;
\item"(ii)" for two $z,z'$ in $M$, not in $2M$, 
with 
$\langle z \mid z \rangle \equiv 0 \pmod 4$
and $\langle z' \mid z' \rangle \equiv 0 \pmod 4$,
$M^z = M^{z'}$ if and only if $z'-z \in 2M$
and 
$\langle z \mid z \rangle \equiv 
\langle z' \mid z' \rangle \pmod 8$;
\item"(iii)" any integral unimodular lattice $M'$
in $V$ such that $M \cap M'$ is of index $2$ in both $M$ and $M'$
appears as one of the lattices $M^z$
($z \in M$, $z \notin 2M$, $\langle z \mid z \rangle \equiv 0 \pmod 4$).
\endroster
(For the analogous facts concerning an odd lattice $M$,
see \cite{Venkov--79}).
% 
% For all this, see e.g. \cite{??????????????????????????????}
%
\enddemo

\demo{Step three: short representatives in $M$ of non-zero
isotropic classes in $\overline M$}

   Let $V$ be again of dimension $8$,
so that two even unimodular lattices in $V$ are always isomorphic,
and let $M$ be such a lattice. 
The $2^8$ elements of $\overline M$ splits as
\roster
\item"(a)" the origin,
\item"(b)" $120$ elements represented by pairs $\pm r$ of roots,
   which are nonisotropic for $q$,
\item"(c)" $135$ nonzero isotropic elements.
\endroster
Let $\psi : M_4 \longrightarrow \overline M$ be the restriction to $M_4$
of the canonical projection. 
We claim that each fiber of $\psi$ has at most $16$ elements. 
Indeed, for $x_1,x_2 \in M_4$ such that $x_1 \ne x_2$ and
$\overline{x_1} = \overline{x_2}$, 
there exists $m \in M$ such that $x_2-x_1=2m \ne 0$ 
and, upon changing $x_2$ to $-x_2$ if necessary, 
$\langle x_1 \mid x_2 \rangle \ge 0$. Then
\roster
\item"(i)" 
$\langle x_2-x_1 \mid x_2-x_1 \rangle = 
4+4-2\langle x_1 \mid x_2 \rangle \le 8$,
\item"(ii)"
$\langle x_2-x_1 \mid x_2-x_1 \rangle = 
4\langle m \mid m \rangle \ge 8$,
\endroster
so that $\langle x_2-x_1 \mid x_1-x_2 \rangle = 8$ and $x_1 \perp x_2$.
The claim follows.

Since $M_4$ has exactly $2160 = 135 \times 16$ elements and since
the image of $\psi$ is {\it a priori} inside the set
of $135$ nonzero isotropic elements of $\overline M$,
it follows that $\psi$ is onto this set.
\enddemo

\demo{Step four: Let $M,M' \subset V \approx \Bbb R^8$ be two even
unimodular lattices such that $M \cap M'$ is of index $2$ in both $M$ and
$M'$. Then there exists $x \in M_4$ such that $M' = s_x(M)$}

Choose a nonzero isotropic element
$\overline z \in \overline M$
represented by an element $x \in M_4$ (see Step three).
Since $M$ is unimodular, there exists a root $r \in M_2$
such that $\langle r \mid x \rangle = 1$.
The element $z = x - 2r$ is also a representative of $\overline z$, and
$$
\langle z \mid z \rangle
\, = \,  \langle x \mid x \rangle 
- 4 \langle r \mid x \rangle 
+ 4 \langle r \mid r \rangle
\, = \,  4 - 4 + 8 = 8 .
$$
Thus $M^z$ is the even neighbour $M'$ of $M$ (see Step two).

   To conclude the proof of Step four, 
we have to check that $s_x(M) = M^z$.
For $m \in M$, we have
$$
\aligned
s_x(m) \, &= \, 
m - \frac{1}{2}\langle z+2r \mid m \rangle (z+2r) \\
\, &= \,
m - \frac{1}{2} \langle z \mid m \rangle z
- \langle z \mid m \rangle r 
- \langle r \mid m \rangle z
- 2 \langle r \mid m \rangle r
\endaligned
$$
and
$$
\langle z \mid s_x(m) \rangle \, \equiv \,
\langle z \mid m \rangle
- \langle z \mid m \rangle \langle z \mid r \rangle
\pmod{2} .
$$
If $\langle z \mid m \rangle$ is even, 
namely if $m \in M_z$,
these formulas show that $s_x(m) \in M_z$,
and consequently that $s_x(m) \in M^z$.
If $\langle z \mid m \rangle$ is odd, 
namely if $m \in M \setminus M_z$,
they show that $s_x(m) - \frac{1}{2}z \in M_z$,
and we have again $s_x(m) \in M^z$.
Thus $s_x(M) \subset M^z$.
Since $s_x(M)$ and $M^z$ are both unimodular,
this shows that $s_x(M) = M^z$.

[Conversely, for any $x \in M_4$, we have 
$\delta(M,s_x(M)) = 1$.]
\enddemo 

\demo{End of proof of Proposition~5.1}

Consider, as in the beginning of this proof, 
an element $g \in O(8,\Bbb Z[1/2])$,
the lattice $M = g(L) \in \Cal L$ at distance $d$ from $L$,
and a sequence $\left( M^{(j)} \right)_{0 \le j \le d}$
such that $M^{(0)} = L$, $M^{(d)} = M$, and
$M^{(j-1)} \cap M^{(j)}$ of index $2$  
in both $M^{(j-1)}$ and $M^{(j)}$.
The previous steps show that there exists
for each $j \in \{1,\hdots,d\}$ an element
$x_j \in \left(M^{(j-1)}\right)_4$ such that
$s_{x_j}(M^{(j-1)}) = M^{(j)}$.
Set $g_j = \prod_{i=1}^{j-1}s_{x_i}$.
Since $g^{-1}g_d$ is in $\operatorname{Aut}(L)$,
and in particular in $O(8,\Bbb Z[1/2])$,
by Step one, 
it is enough to show that 
$g_d \in O(8,\Bbb Z[1/2])$.
We claim that $g_j  \in O(8,\Bbb Z[1/2])$
for all $j \in \{1,\hdots,d\}$,
and we prove the claim by induction on $j$.

The claim is clear for $j = 1$. Assume it holds for some value
of $j \in \{1,\hdots,d-1\}$, and let us check it for $j+1$.
Let $y_{j} \in L$ be the element such that
$g_j(y_{j}) = x_j \in M^{(j)}$.
Then $g_{j+1} = s_{x_j}g_j = g_js_{y_j}g_j^{-1}g_j$
is indeed in $G$ since $g_j \in G$ and $s_{y_j} \in \Gamma_{L,4}$. 
This shows that the claim holds for $j+1$.
$\square$
\enddemo

\bigskip

\proclaim{5.2.\ Corollary} Let $L$ be a root lattice of
type $E_8$. 
Then the subgroup $SO(8,\Bbb Z[1/2])$
of elements of determinant $+1$ in
$O(8,\Bbb Z[1/2])$
is generated by the products $s_xs_{x'}$,
where $x,x'$ are in the shell $L_4$ of $L$.
\endproclaim

\bigskip

\noindent {\bf 5.3.\ A set of $8$ generators for $O(8,\Bbb Z[1/2])$.}
   Let $(\epsilon_1,\hdots,\epsilon_8)$ denote the canonical orthonormal
basis in $\Bbb R^8$. As in \cite{Bourb--68}, set
$$
\aligned
\alpha_1 \, &= \, \frac{1}{2}(\epsilon_1 + \epsilon_8) -
    \frac{1}{2}(\epsilon_2 + \cdots + \epsilon_7) \\
\alpha_2 \, &= \, \epsilon_1 + \epsilon_2 \\
\alpha_j \, &= \, \epsilon_{j-1} - \epsilon_{j-2}
   \quad (j = 3,\hdots,8)
\endaligned
$$
so that $\{\alpha_1,\hdots,\alpha_8\}$ is a basis of a root system of type
$E_8$, and therefore also a basis of a root lattice $L \subset \Bbb R^8$ of
type $E_8$.
Let $s_1,\hdots,s_8$ denote the reflections associated to
$\alpha_1,\hdots,\alpha_8$.
Then $\operatorname{Aut}(L)$ has a Coxeter presentation with
generators $s_1,\hdots,s_8$ and relations of the familiar form
$(s_is_j)^{m_{i,j}} = 1$.

   Consider the vector $2\epsilon_2 \in L_4$ and the corresponding reflection
$\tilde s_2 : x \longmapsto x - \langle x \mid \epsilon_2 \rangle \epsilon_2$.
A simple computation shows that the conjugation by $\tilde s_2$
exchanges $s_2$ with $s_3$ 
and leaves $s_j$ invariant for $j = 1, 4, 5, 6, 7, 8$.
It follows from Proposition~5.1 that $O(8,\Bbb Z[1/2])$ 
has a generating set obtained from that of $\operatorname{Aut}(L)$
by replacing $s_2$ by $\tilde s_2$.
Moreover, the order $\tilde m_{2,j}$ of $\tilde s_2 s_j$ 
in $O(8,\Bbb Z[1/2])$ is equal to
$$
\aligned
\infty \qquad &\text{for} \qquad j = 1 \\
4 \qquad &\text{for} \qquad j = 3,4 \\
3 \qquad &\text{for} \qquad j = 5, 6, 7, 8. 
\endaligned
$$
Thus,  for the generators of $O(8,\Bbb Z[1/2])$ described here, 
the orders $\tilde m_{i,j}$ of the products of two generators
coincide with the corresponding orders $m_{i,j}$ 
for the Coxeter generators
of $\operatorname{Aut}(L)$, 
this for all but three pairs of indices, 
namely for all but $(2,1)$, $(2,3)$, and $(2,4)$.

\enddocument